\documentclass[journal]{IEEEtran}
\ifCLASSINFOpdf

\else

\fi

\usepackage{cite}
\usepackage{psfrag}
\usepackage{rotating}
\usepackage[usenames]{color}
\usepackage{graphicx}
\graphicspath{{figures/}}
\usepackage{epsfig}    
\usepackage{amssymb}   
\usepackage[cmex10]{amsmath}   
\usepackage[font=footnotesize,caption=false]{subfig}
\usepackage{fixltx2e}
\usepackage{array}
\usepackage{dblfloatfix}
\usepackage[export]{adjustbox}

\usepackage{url}
\usepackage{tikz}

\usepackage{algorithm}     
\usepackage{algpseudocode} 

\DeclareMathOperator*{\minimize}{minimize}

\DeclareMathOperator*{\minimum}{min}
\DeclareMathOperator*{\subject}{subject\ to}

\DeclareMathOperator*{\argmin}{arg\ min}

\DeclareMathOperator*{\diag}{diag}

\DeclareMathOperator*{\vectorize}{vec}

\DeclareMathOperator*{\find}{find}

\DeclareMathOperator*{\IQC}{IQC}
\DeclareMathOperator*{\children}{ch}

\DeclareMathOperator*{\real}{Re}
\DeclareMathOperator*{\imag}{Im}
\DeclareMathOperator*{\svec}{svec}
\DeclareMathOperator*{\smat}{smat}
\DeclareMathOperator*{\mat}{mat}

\DeclareMathOperator*{\blkdiag}{blk\ diag}

\newcounter{thm}

\newtheorem{theorem}[thm]{Theorem}

\newcounter{remcount}
\newtheorem{rem}[remcount]{Remark}

\definecolor{red}{rgb}{1,0,0}

\begin{document}
%
\title{Distributed Semidefinite Programming with Application to Large-scale System Analysis}
%
%
%

\author{Sina~Khoshfetrat~Pakazad,~\IEEEmembership{Member,~IEEE,}
        Anders~Hansson,~\IEEEmembership{Member,~IEEE,}
        Martin~S.~Andersen,~\IEEEmembership{Member,~IEEE,}
        and~Anders~Rantzer,~\IEEEmembership{Member,~IEEE}
\thanks{S. Khoshfetrat Pakazad and  A. Hansson are with the Division of Automatic Control, Department of Electrical Engineering, Link\"oping University, Sweden. Email: \{sina.kh.pa, hansson\}@isy.liu.se.}
\thanks{M. S. Andersen is with Department of Applied Mathematics and Computer Science, Technical University of Denmark, Denmark. Email: mskan@dtu.dk}
\thanks{A. Rantzer is with Department of Automatic Control, Lund University, Sweden. Email: anders.rantzer@control.lth.se}
}

\maketitle

\begin{abstract}
Distributed algorithms for solving coupled semidefinite programs (SDPs) commonly require many iterations to converge. They also put high computational demand on the computational agents. In this paper we show that in case the coupled problem has an inherent tree structure, it is possible to devise an efficient distributed algorithm for solving such problems. This algorithm can potentially enjoy the same efficiency as centralized solvers that exploit sparsity. The proposed algorithm relies on predictor-corrector primal-dual interior-point methods, where we use a message-passing algorithm to compute the search directions distributedly. Message-passing here is closely related to dynamic programming over trees. This allows us to compute the exact search directions in a finite number of steps. Furthermore this number can be computed a priori and only depends on the coupling structure of the problem. We use the proposed algorithm for analyzing robustness of large-scale uncertain systems distributedly. We test the performance of this algorithm using numerical examples.
\end{abstract}

\begin{IEEEkeywords}
SDPs, distributed algorithms, primal-dual methods, robustness analysis, interconnected uncertain systems.
\end{IEEEkeywords}

%
\IEEEpeerreviewmaketitle

\section{Introduction}\label{sec:introduction}

Semidefinite programs are convex optimization problems that include linear matrix inequalities (LMIs) or semidefinite constraints. The computational complexity of solving such problems commonly scales badly with the number of optimization variables and/or the dimension of the semidefinite constraints in the problem. This limits our ability to solve large SDPs. Despite this, large SDPs are appearing more and more in different engineering fields, e.g., in problems related to sensor networks, smart grids and analysis of uncertain systems, e.g., see \cite{bis:06,bai:08,mad:15,and:13,meg:97}. This has been the driving force for devising efficient and tailored centralized solvers for such problems. These solvers exploit the structure in the problem to reduce the computational burden of solving the problem in a centralized manner, see e.g., \cite{kim:09,and:13,van:05,han:00b,wal:09}. Despite the success of such approaches for solving medium to large-scale problems, there are still problems that cannot be solved using centralized solvers, see e.g., \cite{sim:14,and+han:12,bis:08,dal:13}. This can be due to limited available computational power and/or memory that prohibits us from solving the problem. Also it can be due to certain structural constraints, e.g., privacy requirements, that obstructs us from even forming the centralized problem.

For such instances, distributed algorithms may be used for solving the problem. These algorithms facilitate solving the problem using a network of computational agents, without the need for a centralized unit. Due to this, the computational complexity of these algorithms scales better, and they potentially enable us to address structural constraints in the problem. The main approach for designing distributed algorithms consists of two major phases. First the structure in the problem is exploited to decompose the problem or reformulate it as a coupled problem. Then first-order splitting methods are used for solving the resulting problem distributedly, see e.g., \cite{sun:14,lu:07}. This approach has been used in many applications, e.g., see \cite{sim:14,kho:14,dal:13}. In \cite{sim:14} the authors consider a sensor localization problem and use a so-called edge-based decomposition for reformulating the underlying SDP as a coupled one. They then employ alternating direction method of multipliers (ADMM) to solve the problem distributedly. An optimal power flow problem has been considered in  \cite{dal:13}, where the authors reformulate the problem as a coupled SDP using semidefinite relaxation techniques. They then use ADMM to solve the coupled problem distributedly. In \cite{kho:14} the authors consider robustness analysis of large-scale interconnected uncertain systems. They exploit the sparsity in the interconnections to decompose the underlying SDP and reformulate it as a coupled problem. This problem is then solved distributedly using algorithms that rely on proximal splitting methods.

The algorithms designed using the aforementioned approach, although effective, suffer from some issues. For instance, since these algorithms rely on first-order splitting methods, with convergence rates $\mathcal O(1/k)$ or $\mathcal O(1/k^2)$ where $k$ is the number of iterations, they require many iterations to converge to an accurate enough solution. Furthermore, exploiting structure and decomposing problems is commonly done through introduction of consensus constraints, which describe the coupling structure in the problem. The number of such constraints is commonly large for SDPs, which can in turn adversely affect the computational and/or convergence properties. Moreover the agents involved in these distributed algorithms need to solve an SDP at every iteration of the algorithm, which can potentially put a considerable computational burden on the agents.

In this paper we propose a distributed algorithm for solving coupled SDPs with a tree structure. These SDPs are defined in Section~\ref{sec:TreeStructureMP}. This algorithm does not suffer from any of the aforementioned issues. We achieve this by avoiding the use of first-order splitting methods and instead rely on primal-dual interior-point methods, which have superior convergence properties. The proposed algorithm is produced by distributing the computations conducted at each iteration of the primal-dual method. Particularly, we use a message-passing algorithm for computing the search directions. Message passing, here, is closely related to non-serial dynamic programming, \cite{kho:15c,kol:09,ber:73}. We also present a similar approach for distributing the remaining computations at every iteration. As a consequence, at each iteration of the primal-dual method, the computational burden on each agent is very low. In fact during each iteration, an agent is required to factorize a relatively small matrix once and is required to communicate with its neighbors twelve times.

The proposed algorithm in this paper is closely related to that of \cite{kho:15c}. In fact, the authors in \cite{kho:15c} use the same approach for devising a distributed algorithm for solving coupled non-conic problems. However, the computation of search directions for SDPs is not as straightforward as for non-conic problems. This is due to introduction of scaling matrices and their inverses in the KKT system, which destroys the structure in the problem. In order to circumvent this issue, we here put forth a novel way for computing the search directions at each iteration. This in turn enables us to use the message-passing algorithm for computing the search directions.

Notice that by using this approach for computing the search directions, we implicitly solve the so-called augmented system. This is done by computing a block $LDL^T$ factorization of its coefficient matrix using a fixed pivoting ordering, where the ordering is enforced by the coupling structure in the problem, \cite{kho:15c}. This is in contrast to existing methods that commonly solve the so-called Schur complement system or normal equations. As a result, the proposed algorithm provides us with more stable and accurate implementation, \cite{wrig:97,cai:06}. Solving the augmented system is also considered in \cite{mat:12}, where the authors also compute the search directions through solving the augmented system by computing an $LDL^T$ factorization using fixed pivoting ordering. This is particularly done by using regularization and iterative refinement. In this paper, however, a block $LDL^T$ factorization is computed using a fixed pivoting ordering without the use of regularization. Hence, the augmented system is solved without the need for iterative refinement.

We then use the proposed algorithm for analyzing large-scale interconnected uncertain systems, distributedly. This is made possible by exploiting the sparsity in the interconnections, as outlined in \cite{and:13}. A similar approach was also used in \cite{kho:14}. There, the authors utilized the so-called range-space decomposition for reformulating the analysis problem as a coupled feasibility problem. They then used algorithms that rely on proximal splitting methods for solving it distributedly. We here instead use the so-called domain-space decomposition to reformulate the analysis problem as a coupled SDP. The coupling structure of this coupled problem is less complicated than that of in \cite{kho:14}, and has a tree structure. This then enables us to use the presented distributed algorithm for solving the problem efficiently and distributedly. We illustrate the performance of the algorithm using numerical examples.

\subsection*{Outline}

Next we first define some notations that are used throughout the paper. In Section \ref{sec:partial} we put forth a definition of coupled and loosely coupled SDPs. We review a predictor-corrector primal-dual interior-point method in Section \ref{sec:primaldual} and briefly discuss how the structure in coupled problems is reflected in the computations conducted at every iteration of this method. Section \ref{sec:TreeStructureMP} expresses coupled problems with a tree structure and discusses the use of message-passing algorithm for solving coupled problems with a tree structure. This is then used in Section~\ref{sec:Dprimaldual} where we present the proposed distributed algorithm for solving coupled SDPs with tree structure. In Section~\ref{sec:ChordalSDP} we discuss a decomposition approach for sparse SDPs. This approach is used in Section \ref{sec:RSA} for reformulating the problem of robustness analysis of large-scale interconnected uncertain systems as coupled SDPs with a tree structure. We test the performance of the proposed distributed algorithm when applied to this problem using numerical experiments in Section~\ref{sec:numerical}. Finally we finish the paper with some concluding remarks in Section~\ref{sec:conclusions}.

\subsection*{Notation}

We denote the set of real and complex numbers with $\mathbb R$ and $\mathbb C$, and the set of $m \times n$ real and complex matrices with $\mathbb R^{m \times n}$ and $\mathbb C^{m \times n}$, respectively. The transpose and conjugate transpose of a matrix $X$ is denoted by $X^T$ and $X^{\ast}$, respectively. The null space of a matrix $X$ is denoted by $\mathcal N(X)$. With $\mathbb S^n$ and $\mathbb H^n$ we denote the set of $n \times n$ symmetric and Hermitian matrices. The set of integer numbers $\{ 1, \dots, n \}$ is denoted by $\mathbb N_n$. Given a set of positive integers $J \subseteq \mathbb N_n$, the matrix $E_J \in \mathbb R^{|J|\times n}$ is a 0--1 matrix obtained from an $n \times n$ identity matrix with rows indexed by $\mathbb N_n \setminus J$ removed, where $|J|$ denotes the number of elements in $J$. This means that $E_J x$ is a $|J|$-dimensional vector that contains the elements of $x$ indexed by $J$. We denote this vector by $X_{_{J}}$. By $x^{i,(k)}_l$ and $X^{i,(k)}_{mn}$ we denote the $l$th element of vector $x^i$ and the element at row $m$ and column $n$ of matrix $X^i$ at the $k$th iteration, respectively. Given matrices $X^k$ for $k = 1, \dots, N$, $\blkdiag(X^1, \dots, X^N)$ denotes a block-diagonal matrix with blocks specified by the given matrices. Similarly $\diag(x_1, \dots, x_N)$ is a diagonal matrix with diagonal elements $x_1, \dots, x_N$. Given vectors $x^k$ for $k= 1, \dots, N$, the column vector $(x^1, \dots, x^N)$ is all of the given vectors stacked. The generalized matrix inequality $G \prec H$ ($G \preceq H$) means that $G-H$ is negative (semi)definite. Given a matrix $X \in \mathbb R^{m\times n}$, $\vectorize(X)$ is an $mn$-dimensional vector that is obtained by stacking all columns of $X$ on top of each other. Given two matrices $X, Y \in \mathbb R^{m \times n}$, $X \bullet Y := \vectorize(X)^T\vectorize(Y)$. For a symmetric matrix $X \in \mathbb S^n$
\begin{multline*}
\svec(X) := (X_{11}, \sqrt{2}X_{21}, \dots, \sqrt{2} X_{n1}, X_{22},\\ \sqrt{2} X_{32}, \dots, \sqrt{2} X_{n2}, \dots, X_{nn}).
\end{multline*}
Operators $\mat$ and $\smat$ are defined as inverses of $\vectorize$ and $\svec$, respectively. Given two matrices $X$ and $Y$ by $X\otimes Y$ we denote the standard Kronecker product. Given $X \in \mathbb S^n$, define $U$ as an $n(n + 1)/2 \times n^2$ matrix such that $U \vectorize(X) = \svec(X)$. Then for two matrices $X, Y \in \mathbb R^{n\times n}$, $\otimes_s$ denotes the symmetrized Kronecker product that is defined as
\begin{align*}
X \otimes_s Y := \frac{1}{2} U(X \otimes Y + Y \otimes X)U^T.
\end{align*}
For properties of the symmetrized Kronecker product refer to \cite{tod:98}. Given two sets $J_1$ and $J_2$, $J_1 \times J_2$ denotes the standard cartesian product and by $J_1 \times_s J_1$ we denote the symmetrized cartesian product defined as
\begin{align*}
J_1 \times_s J_1 := \{ (j,k) \in J_1 \times J_1 \ | \  j\leq k \}.
\end{align*}
For these two sets $J_1 \setminus J_2$ denotes the standard set minus. By $\minimum$ we denote the minimum value and with $\argmin$ we denote the minimizing argument of a function. By $\mathcal L_2^n$ we denote the set of $n$-dimensional square integrable signals, and $\mathcal{RH}_{\infty}^{m \times n}$ represents the set of real, rational $m \times n$ transfer function matrices with no poles in the closed right half plane. A graph is denoted by $Q(V,\mathcal E)$ where $V = \{v_1, \dots, v_n\}$ is its set of vertices or nodes and $\mathcal E \subseteq V\times V$ denotes its set of edges. An induced graph by $V^\prime \subseteq V$ on $Q(V,\mathcal E)$, is a graph $Q_I(V^\prime,\mathcal E^\prime)$ where $\mathcal E^\prime = \mathcal E\cap (V^\prime \times V^\prime)$.

\section{Coupled and Loosely Coupled SDPs}\label{sec:partial}
Let us consider a coupled SDP given~as
\small
\begin{subequations}\label{eq:coupledSDP}
\begin{align}
\minimize_{X} & \quad \sum_{i = 1}^N W^i \bullet X_{_{J_iJ_i}}\\
\subject & \quad Q_j^i\bullet X_{_{J_iJ_i}} = b_j^i, \quad j = 1, \dots, m_i,\notag\\
& \hspace{32mm} i = 1, \dots, N,\\
& \quad X_{_{J_iJ_i}} \succeq 0, \quad i = 1, \dots, N,
\end{align}
\end{subequations}
\normalsize
where $Q_j^i, W^i \in \mathbb S^{|J_i|}$ such that $\begin{bmatrix} \svec(Q_1^i) & \dots & \svec(Q_{m_i}^i) \end{bmatrix}$ has full column rank for all $i = 1, \dots, N$, with the ordered sets $J_i \subseteq \mathbb N_n$ such that $\bigcup_{i=1}^N J_i = \mathbb N_n$, and $X_{_{J_iJ_i}} = E_{J_i}XE_{J_i}^T$ with $X\in \mathbb S^n$ such that $X_{jk} = 0$ if $(j,k) \notin \mathcal J$ and $\mathcal J = \bigcup_{i=1}^N \mathcal J_i := J_i \times_s J_i$. This problem can be seen as a combination of $N$ coupled subproblems, each of which defined by the objective function $W^i\bullet X_{_{J_iJ_i}}$ and constraints $Q_j^i\bullet X_{_{J_iJ_i}} = b_j^i$ for $j = 1, \dots, m_i$ and $X_{_{J_iJ_i}}\succeq 0$. Let us now define $\mathcal I_{(i,j)} = \left\{ k \ | \ (i,j) \in \mathcal J_k \right\}$, which denotes the set of subproblems that are coupled in that they all depend on the variable $X_{ij}$. Notice that agents $a$ and $b$ are members of $\mathcal I_{(i,j)}$ if and only if $\{ i, j\} \subseteq J_a \cap J_b$. It is possible to provide a more explicit description of the coupling among the subproblems by decomposing \eqref{eq:coupledSDP} as
\small
\begin{subequations}\label{eq:coupledSDPDecomposed}
\begin{align}
\minimize_{X,\bar X^i} & \quad \sum_{i = 1}^N W^i \bullet \bar X^i\label{eq:coupledSDPDecomposed-a}\\
\subject & \quad Q_j^i\bullet \bar X^i = b_j^i, \quad i = 1, \dots, N,\label{eq:coupledSDPDecomposed-b}\\
& \quad \bar X^i \succeq 0, \quad i = 1, \dots, N,\label{eq:coupledSDPDecomposed-c}\\
& \quad  \bar X^i = E_{J_i}XE_{J_i}^T, \quad i = 1, \dots, N\label{eq:coupledSDPDecomposed-d}.
\end{align}
\end{subequations}
\normalsize
Notice that in \eqref{eq:coupledSDPDecomposed}, the objective function terms and constraints in~\eqref{eq:coupledSDPDecomposed-a}--\eqref{eq:coupledSDPDecomposed-c} are decoupled and the coupling in the problem is described using the consensus constraints in \eqref{eq:coupledSDPDecomposed-d}. It is also possible to provide a graphical representation of the coupling using undirected graphs. Particularly let $Q_s(\mathcal J,\mathcal E_s)$ be a graph with vertex set $\mathcal J$ as defined above and edge set $\mathcal E_s = \left\{ \left( (i,j),(v,t) \right) \  | \  \mathcal I_{(i,j)} \cap \mathcal I_{(v,t)} \neq \emptyset   \right\}$. We refer to this graph as the sparsity graph of the problem. Let us now illustrate the definitions above using an example given as
\small
\begin{subequations}\label{eq:example}
\begin{align}
\minimize_{X} & \quad W^1\bullet X_{\{1,2,4\}\{1,2,4\}} + \notag\\  & \quad \quad \quad  W^2\bullet X_{\{1,3,4\}\{1,3,4\}} + W^3\bullet X_{\{4,5\}\{4,5\}}\\
\subject & \quad \begin{bmatrix} x_{11} & x_{12} & x_{13} & 0 & 0\\ x_{12} & x_{22} & 0 & x_{24} & 0 \\ x_{13} & 0 & x_{33} & x_{34} & 0 \\ 0 & x_{24} & x_{34} & x_{44} & x_{45} \\ 0 & 0 & 0 & x_{45} & x_{55} \end{bmatrix} \succeq 0\label{eq:example-b}.
\end{align}
\end{subequations}
\normalsize
Notice that the constraint in \eqref{eq:example-b} can be rewritten as
\small
\begin{multline*}
\begin{bmatrix} x_{11} & x_{12} & x_{13} & 0 & 0\\ x_{12} & x_{22} & 0 & x_{24} & 0 \\ x_{13} & 0 & x_{33} & x_{34} & 0 \\ 0 & x_{24} & x_{34} & x_{44} & x_{45} \\ 0 & 0 & 0 & x_{45} & x_{55} \end{bmatrix} = \\ E_{J_1}^T  \begin{bmatrix} \frac{x_{11}}{2} & x_{12} & 0 \\x_{12} & x_{22} & x_{24} \\ 0 & x_{24} & \frac{x_{44}}{3} \end{bmatrix} E_{J_1}\\ \quad \quad \quad+ E_{J_2}^T  \begin{bmatrix} \frac{x_{11}}{2} & x_{13} & 0 \\x_{13} & x_{33} & x_{34} \\ 0 & x_{34} & \frac{x_{44}}{3} \end{bmatrix} E_{J_2} \\ + E_{J_3}^T  \begin{bmatrix} \frac{x_{44}}{3} & x_{45} \\x_{45} & x_{55} \end{bmatrix} E_{J_3} \succeq 0.
\end{multline*}
\normalsize
with $J_1 = \{ 1, 2, 4\}$, $J_2 = \{ 1, 3, 4\}$ and $J_3 = \{ 4, 5 \}$. Then the optimal objective value of

\begin{figure}[t]
\begin{center}
\includegraphics[width=4.5cm]{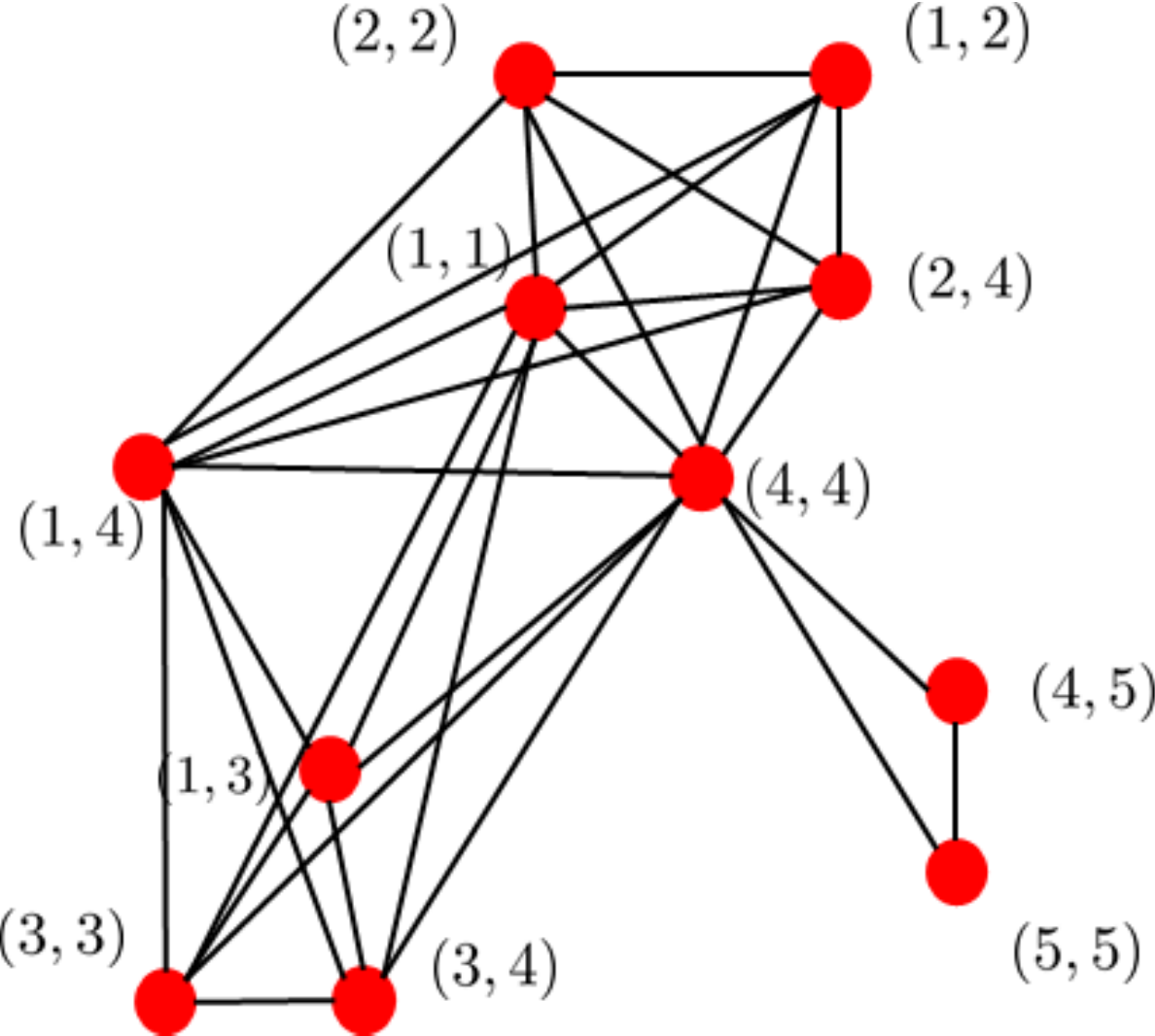}    
\caption{\small Sparsity graph for the coupled SDP in \eqref{eq:exampleCoupled}.\normalsize }
\label{fig:ex1}
\end{center}
\end{figure}

\small
\begin{subequations}\label{eq:exampleCoupled}
\begin{align}
\minimize_{X} & \quad W^1\bullet X_{\{1,2,4\}\{1,2,4\}} + \notag\\  & \quad \quad  \quad W^2\bullet X_{\{1,3,4\}\{1,3,4\}} + W^3\bullet X_{\{4,5\}\{4,5\}}\\
\subject & \quad \begin{bmatrix} \frac{x_{11}}{2} & x_{12} & x_{14} \\x_{12} & x_{22} & x_{24} \\ x_{14} & x_{24} & \frac{x_{44}}{3} \end{bmatrix} \succeq 0, \quad x_{14} = 0,\\
& \quad \begin{bmatrix} \frac{x_{11}}{2} & x_{13} & x_{14} \\x_{13} & x_{33} & x_{34} \\ x_{14} & x_{34} & \frac{x_{44}}{3} \end{bmatrix} \succeq 0, \quad x_{14} = 0,\\
& \quad \begin{bmatrix} \frac{x_{44}}{3} & x_{45} \\x_{45} & x_{55} \end{bmatrix}\succeq 0,
\end{align}
\end{subequations}
\normalsize
defines an upperbound for the optimal objective value of \eqref{eq:example}. The problem in \eqref{eq:exampleCoupled} is a coupled SDP with smaller semidefinite constraints. This method for reformulating the problem is commonly used for cases when the original problem is either impossible or very difficult to solve. Notice that this problem is in the same format as~\eqref{eq:coupledSDP}. The sparsity graph of this problem is illustrated in Figure \ref{fig:ex1}, where for instance there is an edge between the nodes $(1,1)$ and $(1,2)$ since the intersection between the sets $\mathcal I_{(1,1)} = \{ 1, 2 \}$ and $\mathcal I_{(1,2)} = \{ 2 \}$ is nonempty.

\noindent In case for a coupled problem
\begin{itemize}
\item $|J_i \cap J_j|\ll n$ for all $i, j \in \mathbb N_N$;
\item $|\mathcal I_{(i,j)} \cap \mathcal I_{(v,t)}| \ll N$ for all $(i,j), (v,t) \in  \mathcal J$,
\end{itemize}
then we call this problem loosely coupled. As we will see later, it is possible to devise efficient distributed solvers based on primal-dual interior-point methods for solving coupled and loosely coupled SDPs. To this end, let us first briefly review primal-dual interior-point methods for solving SDPs.

\section{Primal-Dual Interior-point Methods for Solving SDPs}\label{sec:primaldual}
It is possible to iteratively solve a standard-form SDP, given~as
\small
\begin{equation}\label{eq:SDP}
\begin{split}
\minimize_{X}  & \quad C \bullet X \\
\subject   & \quad A_i \bullet X = b_i, \quad i = 1, \dots, m,\\
& \quad X \succeq 0,
\end{split}
\end{equation}
\normalsize
where $b \in \mathbb R^m$ and $X, A_i, C \in \mathbb S^n$ such that $\begin{bmatrix} \svec(A_1) & \dots & \svec(A_m) \end{bmatrix}$ has full column rank, using primal-dual interior-point methods. Particularly, given the iterates $(X^{(k)} \succ 0, S^{(k)} \succ 0, v^{(k)})$, a primal-dual interior-point method generates the next iterates $(X^{(k+1)}, S^{(k+1)}, v^{(k+1)})$ by taking a single Newton step applied to the perturbed KKT conditions
\small
\begin{subequations}\label{eq:SDPOptimality}
\begin{align}
A_i \bullet X &= b_i, \quad i = 1, \dots, m,\\
\sum_{i=1}^m v_i A_i + S &= C,\\
XS &= \delta I,\label{eq:SDPOptimality-c}
\end{align}
\end{subequations}
\normalsize
together with $S \succ 0$ and $X \succ 0$ where $\delta > 0$. Specifically this Newton step can be computed by solving the following linear system of equations
\small
\begin{subequations}\label{eq:SDPOptimalityLinearized}
\begin{align}
A_i \bullet \Delta X &= b_i - A_i \bullet X^{(k)}, \ \ i=1, \dots, m,\\
\sum_{i=1}^m \Delta v_i A_i + \Delta S &= C - S^{(k)} - \sum_{i=1}^m v_i^{(k)}A_i, \\
H_D(\Delta X S^{(k)} + X^{(k)} \Delta S) &= \delta I - H_D(X^{(k)}S^{(k)}), \label{eq:SDPOptimalityLinearized-c}
\end{align}
\end{subequations}
\normalsize
where $H_D(M) = 1/2(DMD^{-1} + D^{-T}MD^{T})$, $\delta = \sigma \mu$ is the perturbation parameter with $\mu = X^{(k)}\bullet S^{(k)}/n$ denoting the surrogate duality gap and $\sigma \in [0, 1]$, and where \eqref{eq:SDPOptimalityLinearized-c} is a modified linearization of \eqref{eq:SDPOptimality-c} that ensures that the computed directions $\Delta S$ and $\Delta X$ are symmetric. There are different choices for the scaling matrix $D$ in \eqref{eq:SDPOptimalityLinearized-c}, e.g., see \cite{tod:98} and references therein. For the sake of brevity, we limit our discussion to the choices presented in \cite{nes:97,nes:95}, that is we choose $D = G^{-1}$ with $W = GG^T$ where
\small
\begin{equation}\label{eq:scaling}
\begin{split}
W :&= (X^{(k)})^{\frac{1}{2}} \left( (X^{(k)})^{\frac{1}{2}} S^{(k)} (X^{(k)})^{\frac{1}{2}} \right)^{-\frac{1}{2}} (X^{(k)})^{\frac{1}{2}}\\
& = (S^{(k)})^{-\frac{1}{2}} \left( (S^{(k)})^{\frac{1}{2}} X^{(k)} (S^{(k)})^{\frac{1}{2}} \right)^{\frac{1}{2}} (S^{(k)})^{-\frac{1}{2}}.
\end{split}
\end{equation}
\normalsize
This scaling is referred to as the Nesterov-Todd or NT scaling.  In order to make the notation less complicated, from now on we drop the iteration index $k$, and we use lowercase notation for denoting vectorized variables or residuals, e.g., we use $\Delta x$ as $\svec(\Delta X)$ or $r_{\textrm{dual}}$ as $\svec(R_{\textrm{dual}})$. Using symmetrized Kronecker product we can then rewrite \eqref{eq:SDPOptimalityLinearized} more compactly as
\begin{align}\label{eq:SDPOptimalityLinearizedCompact}
\begin{bmatrix} 0 & A & 0 \\ A^T & 0 & I \\ 0 & U & F \end{bmatrix}\begin{bmatrix} \Delta v \\ \Delta x \\ \Delta s \end{bmatrix} = \begin{bmatrix} r_{\textrm{primal}} \\ r_{\textrm{dual}} \\ r_{\textrm{cent}} \end{bmatrix},
\end{align}
where $A = \begin{bmatrix} \svec(A_1) & \dots & \svec(A_m) \end{bmatrix}^T$, $U = D \otimes_s D^{-T}S$, $F = D X \otimes_s D^{-T}$ and
\begin{equation}
\begin{split}
r_{\textrm{primal}} &= b - Ax\\
R_{\textrm{dual}} &= C - S - \sum_{i=1}^m v_iA_i,\\ R_{\textrm{cent}} &= \delta I - H_D(XS),
\end{split}
\end{equation}
see \cite{tod:98}. One way of solving \eqref{eq:SDPOptimalityLinearizedCompact}, is to first solve for $\Delta s$ as in
\begin{align}\label{eq:DeltaS}
\Delta s = F^{-1}\left( r_{\textrm{cent}}-U \Delta x \right),
\end{align}
and then solve
\begin{align}\label{eq:SDPOptimalityLinearizedCompactAug}
\begin{bmatrix} -F^{-1}U & A^T \\ A & 0 \end{bmatrix} \begin{bmatrix} \Delta x \\ \Delta v  \end{bmatrix} = \begin{bmatrix} r \\ r_{\textrm{primal}}\end{bmatrix}
\end{align}
for $\Delta X$ and $\Delta v$, where $r = r_{\textrm{dual}} - F^{-1}r_{\textrm{cent}}$. Notice that since $F^{-1}U$ is positive definite, \cite[Thm. 3.2]{tod:98},~\eqref{eq:SDPOptimalityLinearizedCompactAug} also describes the optimality condition for the following convex optimization problem
\begin{equation}\label{eq:SDPOptimalityLinearizedCompactAugQP}
\begin{split}
\minimize_{\Delta x} & \quad \frac{1}{2}\Delta x^T F^{-1}U \Delta x + r^T \Delta x \\
\subject & \quad A \Delta x = r_{\textrm{primal}}.
\end{split}
\end{equation}
So it is possible to compute $\Delta X$ and $\Delta v$ by either solving the system of equations in \eqref{eq:SDPOptimalityLinearizedCompactAug} or the problem in~\eqref{eq:SDPOptimalityLinearizedCompactAugQP}. In this paper we focus on predictor-corrector primal-dual methods that rely on modified Newton directions. In order to compute these directions, at each iteration, we need to solve \eqref{eq:SDPOptimalityLinearizedCompactAug} or \eqref{eq:SDPOptimalityLinearizedCompactAugQP} twice with different choices of $r$. We lay out a predictor-corrector primal-dual interior-point method in Algorithm \ref{alg:PD}, based on the work in \cite{tod:98}.
\begin{rem}\label{rem:rem1}
Algorithm \ref{alg:PD} can detect infeasibility of the problem if either the primal or dual iterates diverge. This means that this algorithm is unable to detect weak infeasibility, \cite{dek:00,lou:96}, and generally in such cases converges to a near-feasible solution, \cite{tut:01}.
\end{rem}
The major computational burden of each iteration of this primal-dual method concerns the computation of the predictor and corrector directions. Next we will investigate how the structure in coupled problems is reflected in \eqref{eq:SDPOptimalityLinearizedCompactAugQP} and how this structure can be used to our advantage.
\begin{algorithm}[t]
\caption{Predictor-corrector Primal-dual Interior-point Method, \cite{tod:98}}\label{alg:PD}
\begin{algorithmic}[1]
\small
\State{Given $k = 0$, $\tau \in (0, 1)$, $a \in \{ 1, 2, 3 \}$, $\epsilon>0$, $\epsilon_{\text{feas}}>0$,  initial iterates $(X^{(0)}, S^{(0)}, v^{(0)})$ such that $X^{(0)} \succ 0$ and $S^{(0)} \succ 0$ and $\mu = X^{(0)} \bullet S^{(0)}/ n$.}
\Repeat
\State{Compute $D$.}
\State{Predictor step: Set $\sigma = 0$ and compute the search directions $\Delta X_{\textrm{pred}}$, $\Delta v_\textrm{pred}$ by either solving \eqref{eq:SDPOptimalityLinearizedCompactAug} or \eqref{eq:SDPOptimalityLinearizedCompactAugQP} and $\Delta S_\textrm{pred}$ using~\eqref{eq:DeltaS}.}
\State{Compute primal and dual step sizes as
\begin{equation*}
\begin{split}
\alpha_p &:= \minimum \left( 1, \frac{-\tau}{\lambda_{\minimum}\left((X^{(k)})^{-1}\Delta X_\textrm{pred}\right)} \right),\\
\alpha_d &:= \minimum \left( 1, \frac{-\tau}{\lambda_{\minimum}\left((S^{(k)})^{-1}\Delta S_\textrm{pred}\right)} \right).
\end{split}
\end{equation*}
}
\State{Set $\sigma = \left( \frac{(X^{(k)} + \alpha_p \Delta X_\textrm{pred})\bullet(S^{(k)} + \alpha_d \Delta S_\textrm{pred})}{X^{(k)}\bullet S^{(k)}} \right)^a$.}
\State{Corrector step: Having computed $\sigma$ compute the search directions $\Delta X_\textrm{corr}$, $\Delta v_\textrm{corr}$ by either solving \eqref{eq:SDPOptimalityLinearizedCompactAug} or \eqref{eq:SDPOptimalityLinearizedCompactAugQP} with
\begin{multline*}
r^{(k)} = r_{\textrm{dual}}^{(k)} - (F^{(k)})^{-1}r_{\textrm{cent}}^{(k)} + \\ (F^{(k)})^{-1}\svec(H_D(\Delta X_\textrm{pred}\Delta S_\textrm{pred})),
\end{multline*}
and $\Delta S_\textrm{corr}$ using
\begin{multline*}
\Delta s_\textrm{corr} = (F^{(k)})^{-1}\left( r_{\textrm{cent}}^{(k)}- \right.\\ \left. \svec(H_D(\Delta X_\textrm{pred}\Delta S_\textrm{pred})) -U^{(k)} \Delta x_\textrm{corr} \right).
\end{multline*}
}
\State{Compute primal and dual step sizes as above, though using $\Delta X_\textrm{corr}$ and $\Delta S_\textrm{corr}$.}
\State{Update
\begin{align*}
X^{(k+1)} &=  X^{(k)} + \alpha_p\Delta X_\textrm{corr},\\
S^{(k+1)} &=  S^{(k)} + \alpha_d\Delta S_\textrm{corr},\\
v^{(k+1)} &=  v^{(k)} + \alpha_d\Delta v_\textrm{corr}.
\end{align*}
}
\State {Set $k = k + 1$.}
\State{$\mu = X^{(k)} \bullet S^{(k)}/ n$.}
\Until{$\left\| r_{\textrm{primal}}^{(k)} \right\|^2, \left\| \svec\left(R_{\textrm{dual}}^{(k)}\right) \right\|^2 \leq \epsilon_{\text{feas}}$ and $\mu \leq \epsilon$.}
\normalsize
\end{algorithmic}
\end{algorithm}

Let us apply the primal-dual method in Algorithm \ref{alg:PD} to the coupled SDP in \eqref{eq:coupledSDPDecomposed}. The perturbed KKT optimality conditions for this problem can be written as
\begin{subequations}\label{eq:coupledSDPKKT}
\begin{align}
Q^i_j \bullet \bar X^i &= b^i_j, \ \ j = 1, \dots, m_i,\\
\sum_{j = 1}^{m_i} v^i_j Q^i_j - \smat( \bar v^i) + S^i &= W^i, \\ 
\bar X^i S^i &= \delta I, \\
\bar X^i - X_{_{J_iJ_i}} &= 0,
\end{align}
\end{subequations}
for $i = 1, \dots, N$, together with
\begin{align}
\sum_{i = 1}^N(E_{J_i} \otimes_s E_{J_i})^T \bar v^i &= 0,
\end{align}
and $\bar X^i, S^i\succ 0$ for $i = 1, \dots, N$. Define $\mathcal Q^i = \begin{bmatrix} \svec(Q^i_1) & \dots & \svec(Q^i_{m_i}) \end{bmatrix}^T$. Similar to \eqref{eq:SDPOptimalityLinearized}, given iterates $X$ and $\bar X^{i} \succ 0$ such that they satisfy \eqref{eq:coupledSDPDecomposed-d}, $S^{i}\succ 0$, $v^{i}$ and $\bar v^{i}$ such that $\sum_{i=1}^N(E_{J_i} \otimes_s E_{J_i})^T \bar v^{i} = 0$ for all $i = 1, \dots, N$, the Newton step corresponding to the above system of equations can be computed by solving
\begin{subequations}\label{eq:coupledSDPKKTLinearized}
\begin{align}
&\mathcal Q^i \Delta \bar x^i = b^i - \mathcal Q^i \bar x^{i},\\
&\sum_{j = 1}^{m_i} \Delta v^i_j Q^i_j - \smat( \Delta \bar v^i ) + \Delta S^i = \notag \\  &\hspace{25mm} W^i -  \sum_{j = 1}^{m_i} v^{i}_j Q^i_j + \smat(\bar v^{i}) - S^{i},\\
&H_{D^i}(\Delta \bar X^i S^{i} + \bar X^{i} \Delta S^i) = \delta I - H_{D^i}(\bar X^{i}S^{i}), \\
&\Delta \bar X^i - \Delta X_{_{J_iJ_i}} = 0,
\end{align}
\end{subequations}
for $i = 1, \dots, N$, together with $\sum_{i=1}^N(E_{J_i} \otimes_s E_{J_i})^T \Delta \bar v^i = 0,$ where the scaling matrices $D^i$ are computed as discussed above and in \eqref{eq:scaling}, though based on the given local iterates $\bar X^{i,(k)}$ and $S^{i,(k)}$. This system of equations can be rewritten in a more compact manner as
\begin{align}
\begin{bmatrix}  0 & 0 & \mathbf{\mathcal Q} & 0 & 0 \\ 0 & 0 & I & -\bar {\mathcal  E} & 0\\ \mathbf{\mathcal Q}^T & I & 0 & 0 & I \\ 0 & -\bar{\mathcal E}^T & 0 & 0 & 0 \\ 0 & 0 & \mathbf U & 0 & \mathbf F \end{bmatrix} \begin{bmatrix} \Delta \mathbf v\\ \Delta \bar{\mathbf v} \\ \Delta \bar{\mathbf x} \\ \Delta x \\ \Delta \mathbf s  \end{bmatrix} = \begin{bmatrix} \mathbf r_{\textrm{primal}} \\ 0 \\ \mathbf r_{\textrm{dual}} \\ 0 \\ \mathbf r_{\textrm{cent}} \end{bmatrix}
\end{align}
where $\mathcal Q$, $\mathbf U$ and $\mathbf F$ are block-diagonal with diagonal blocks $\mathcal Q^i$ and
\begin{align*}
U^{i} & = D^i \otimes_s (D^i)^{-T}S^{i}, \\
F^{i} & = D^i \bar X^{i} \otimes_s (D^i)^{-T},
\end{align*}
and $\bar{\mathcal E}^T = \begin{bmatrix} (E_{J_1} \otimes_s E_{J_1})^T & \dots & (E_{J_N} \otimes_s E_{J_N})^T \end{bmatrix}$. Also here $\Delta \mathbf v$, $\Delta \bar{\mathbf v}$, $\Delta \bar{\mathbf x}$ and $\Delta \mathbf s$ denote all the corresponding variables stacked, e.g., $\Delta \mathbf v = \left(\Delta v^1, \dots, \Delta v^N\right)$. Similarly $\mathbf r_{\textrm{primal}}$, $\mathbf r_{\textrm{dual}}$ and $\mathbf r_{\textrm{cent}}$ denote all the primal, dual and centering residuals stacked, where each of the stacked terms in the residual vectors are based on
\begin{subequations}\label{eq:residuals}
\begin{align}
r_{\textrm{primal}}^{i} &= b^i - \mathcal Q^i x^{i},\label{eq:residuals-a}\\
R_{\textrm{dual}}^{i} & = W^i -  \sum_{j = 1}^{m_i} v^{i}_j Q^i_j + \smat(\bar v^{i}) - S^{i}, \label{eq:residuals-b}\\
R_{\textrm{cent}}^{i} & = \delta I - H_{D^i}(\bar X^{i}S^{i}).\label{eq:residuals-c}
\end{align}
\end{subequations}
Similar to before, we compute the primal-dual directions by first solving for $\Delta \mathbf s$ as
\begin{align}\label{eq:DeltaSi}
\Delta s^i = (F^{i})^{-1}\left( r_{\textrm{cent}}^{i} - U^{i}\Delta \bar x^i \right),
\end{align}
or equivalently as
\begin{align}
\Delta s^i = (F^{i})^{-1}\left( r_{\textrm{cent}}^{i} - U^{i} \svec(E_{J_i}^T\Delta XE_{J_i}) \right),
\end{align}
for $i = 1, \dots, N$. Then we solve
\small
\begin{multline}\label{eq:CoupledSDPAugmentedSystem}
\begin{bmatrix} -\mathbf F^{-1}\mathbf U & 0 & \mathcal Q^T & I \\ 0 & 0 & 0 & -\bar{\mathcal E}^T \\ \mathcal Q & 0 & 0 & 0\\ I & -\bar{\mathcal E} & 0 & 0   \end{bmatrix}\begin{bmatrix} \Delta \bar{\mathbf x} \\ \Delta x \\ \Delta \mathbf v\\ \Delta \bar{\mathbf v}\end{bmatrix} = \begin{bmatrix} \mathbf r \\ 0 \\ \mathbf r_{\textrm{primal}} \\ 0 \end{bmatrix},
\end{multline}
\normalsize
where $ \mathbf r = (r^{1}, \dots, r^{N})$ with $r^{i} = r_{\textrm{dual}}^{i} - (F^{i})^{-1}r_{\textrm{cent}}^{i}$. Notice that the system of equations in~\eqref{eq:CoupledSDPAugmentedSystem} also describes the optimality conditions for
\begin{subequations}\label{eq:CoupledSDPQP}
\begin{align}
\minimize_{\Delta \bar{\mathbf x},\Delta x} & \ \ \sum_{i=1}^N\frac{1}{2} (\Delta \bar x^i)^T (F^{i})^{-1}U^{i}\Delta \bar x^i + (r^{i})^T \Delta \bar x^i\\
\subject & \ \ \mathcal Q^i \Delta \bar x^i = r_{\textrm{primal}}^{i}, \quad i = 1, \dots, N,\\
& \ \ \Delta \bar X^i - \Delta X_{_{J_iJ_i}} = 0, \quad i = 1, \dots, N.\label{eq:CoupledSDPQP-c}
\end{align}
\end{subequations}
where $(F^{i})^{-1}U^{i} \succ 0$ for $i = 1, \dots, N$. So the predictor and corrector directions can also be computed by solving~\eqref{eq:CoupledSDPQP}. To be more precise, for the predictor directions, we solve \eqref{eq:CoupledSDPQP}, with $\sigma = 0$, for $\Delta \bar{\mathbf x}_{\textrm{pred}}, \Delta x_{\textrm{pred}}, \Delta \mathbf v_{\textrm{pred}}$ and $\Delta \bar{\mathbf v}_{\textrm{pred}}$, and compute $\Delta \mathbf s_{\textrm{pred}}$ using \eqref{eq:DeltaSi}. For the corrector directions, using the updated $\sigma$, we compute the directions $\Delta \bar{\mathbf x}_{\textrm{corr}}, \Delta x_{\textrm{corr}}, \Delta \mathbf v_{\textrm{corr}}$ and $\Delta \bar{\mathbf v}_{\textrm{corr}}$ by solving \eqref{eq:CoupledSDPQP} with
\begin{multline}\label{eq:RCorrector}
r^{i} = r_{\textrm{dual}}^{i} \\- (F^{i})^{-1}\left(r_{\textrm{cent}}^{i} - \svec(H_{D^i}(\Delta \bar{X}^i_\textrm{pred}\Delta S^i_\textrm{pred})\right)
\end{multline}
and compute $\Delta \mathbf s_{\textrm{corr}}$ as
\begin{multline}\label{eq:DeltaSCorri}
\Delta s_\textrm{corr}^{i} = (F^{i})^{-1}\left( r_{\textrm{cent}}^{i}- \right.\\ \left. \svec(H_{D^i}(\Delta \bar X^{i}_\textrm{pred}\Delta S^i_\textrm{pred})) -U^{i} \Delta \bar x^i_\textrm{corr} \right),
\end{multline}
for $i = 1, \dots, N$. As a result having computed predictor or corrector versions of the directions $\Delta \bar{\mathbf x}, \Delta x, \Delta \mathbf v$ and $\Delta \bar{\mathbf v}$, computing $\Delta s^i_{\textrm{pred}}$ and $\Delta s^i_{\textrm{corr}}$ can be done independently by $N$ computing agents in parallel. Also notice that the coupling structure in \eqref{eq:CoupledSDPQP} is the same as in \eqref{eq:coupledSDPDecomposed}. This allows us to employ distributed computational algorithms to distributedly solve for the search directions using $N$ collaborating agents. To illustrate this, note that the problem in \eqref{eq:CoupledSDPQP} can be written as
\begin{equation}
\begin{split}
\minimize_{\bar x, x} & \quad F_1(\bar x)\\
\subject & \quad A \bar x + B x = c,
\end{split}
\end{equation}
with $\bar x = (\Delta \bar x^1, \dots, \Delta \bar x^N)$ and $x = \Delta x$. This problem can be solved distributedly using proximal splitting methods, e.g., ADMM, \cite{boyd:11,boy:14,ber:97}. The use of proximal splitting methods for computing the primal-dual directions has been considered in \cite{kho:15b,kho:14b}. Devising distributed algorithms for solving coupled SDPs that also rely on this approach can be seen as an extension of the use of the algorithm proposed in \cite{kho:15b} to SDPs. Even though distributed algorithms based on proximal splitting are effective for non-conic problems, they suffer from certain issues when used for solving SDPs. Particularly, notice that the computed search directions using this approach are inexact and first-order splitting methods generally require many iterations to compute accurate enough search directions. Furthermore, the number of consensus constraints in \eqref{eq:CoupledSDPQP-c} are generally large for coupled SDPs which can in turn adversely affect the performance and numerical properties of such splitting methods. Also notice that for a predictor-corrector primal-dual method the search directions are computed through solving a system of the form \eqref{eq:CoupledSDPQP} twice. This means that the iterative scheme for solving \eqref{eq:CoupledSDPQP} needs to be run twice at each iteration of the primal-dual method. Hence, distributed algorithms that rely on proximal or first-order splitting for computing the search directions, potentially, require many iterations to converge to the solution. Despite all such issues, in many cases such splitting methods are among the only resorts for distributedly solving coupled or loosely coupled SDPs. However for coupled problems that have an inherent tree structure, which is common in loosely coupled SDPs, we can  devise an efficient algorithm for solving coupled SDPs. This is the focus of the upcoming sections. But first we express what we mean by the tree structure.

\section{Tree Structure in Coupled Problems and Message Passing}\label{sec:TreeStructureMP}
\begin{figure}[t]
\begin{center}
\includegraphics[width=6cm]{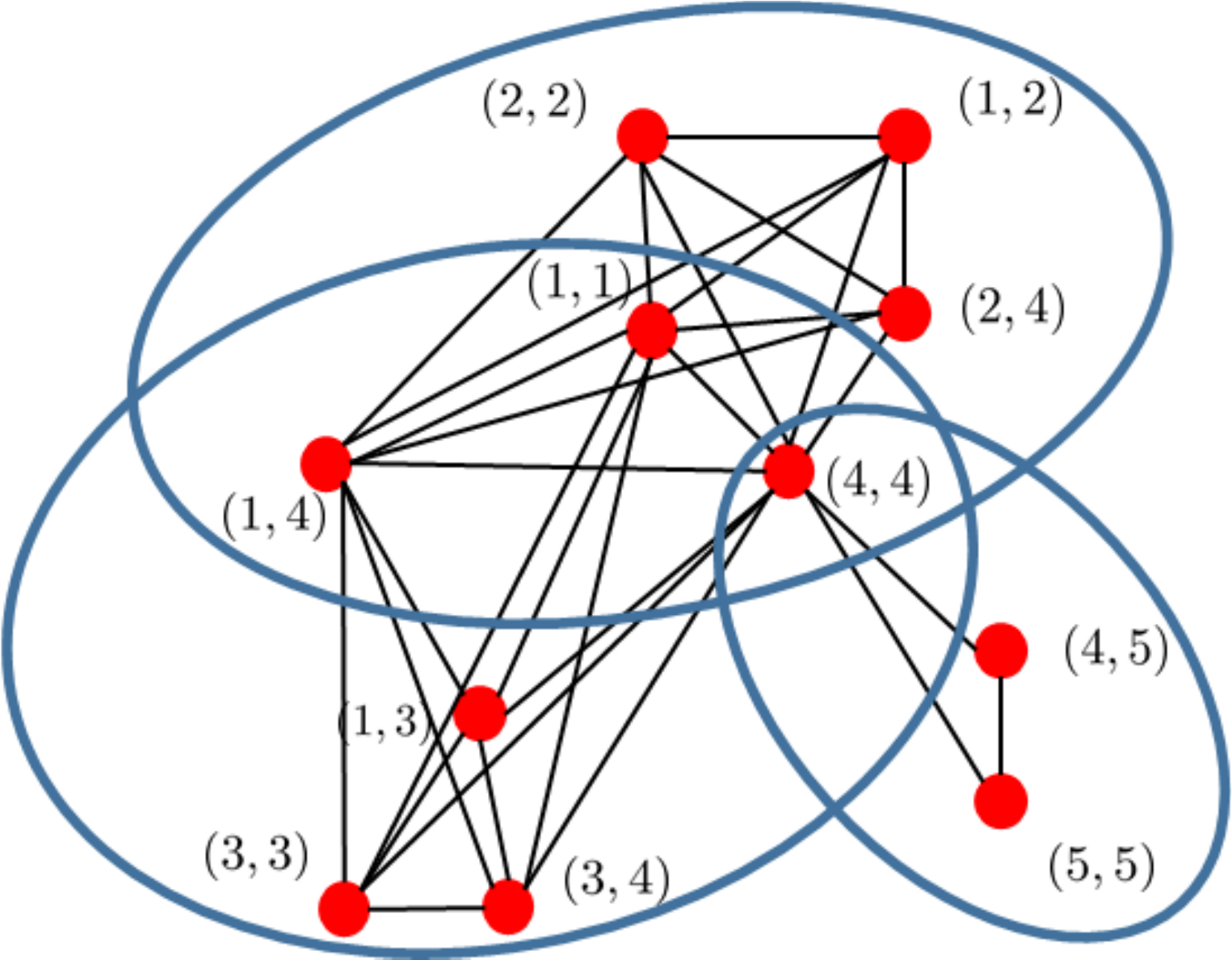}    
\caption{\small Clustered sparsity graph.\normalsize }
\label{fig:ex1b}
\end{center}
\end{figure}
\begin{figure}[t]
\begin{center}
\includegraphics[width=7cm]{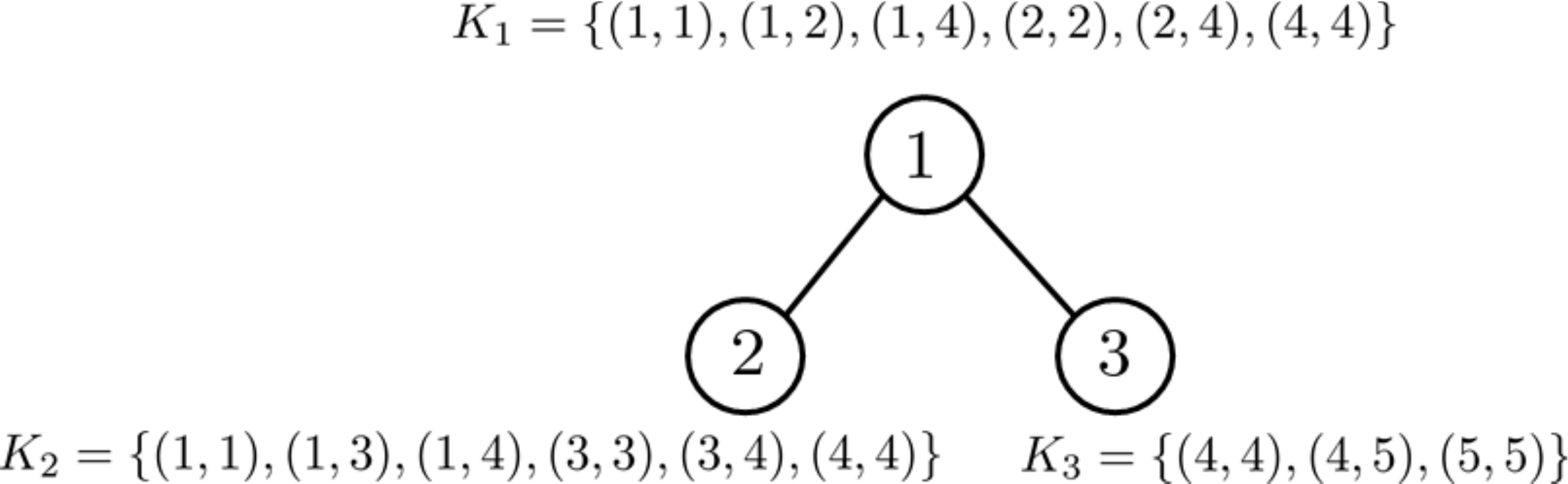}    
\caption{\small Tree representation of the sparsity graph .\normalsize }
\label{fig:ex1c}
\end{center}
\end{figure}
Let us reconsider the coupled SDP in \eqref{eq:exampleCoupled}. Notice that for this problem it is possible to cluster the variables or the nodes in its sparsity graph as shown in Figure~\ref{fig:ex1b}. As can be seen from the figure, each of the clusters induce a complete subgraph on the sparsity graph. We can then provide a more compact representation of the sparsity graph using the tree in Figure~\ref{fig:ex1c}. Each node in this tree corresponds to each of the clusters of variables denoted by $K_i$. Furthermore, for this problem, the tree is such that for every two nodes $i$ and $j$ in the tree, $K_i \cap K_j$ is contained in all the clusters in the path connecting the two nodes in the tree. We refer to problems that enjoy this inherent structure as coupled with a tree structure. Next we lay out an approach for exploiting this structure in coupled problems.

Let us start by describing some definitions relating to graphs. Consider a graph $Q(V,\mathcal E)$. A clique $C_i$ of this graph is a maximal subset of $V$ that induces a complete subgraph on $Q$, i.e., no clique is properly contained in another clique, \cite{blp:94}. Assume that all cycles of length at least four of $Q(V,\mathcal E)$ have a chord, where a chord is an edge between two non-consecutive vertices in a cycle. This graph is then called chordal \cite[Ch.~4]{gol:04}. It is possible to make a non-chordal graph chordal by adding edges to the graph. The resulting graph is then referred to as a chordal embedding. Let $\mathbf C_Q = \{ C_1, \dots, C_q \}$ denote the set of its cliques, where $q$ is the number of cliques of the graph. Then there exists a tree defined on $\mathbf C_Q$ such that for every $C_i, C_j \in\mathbf C_Q$ where $i \neq j$, $C_i \cap C_j$ is contained in all the cliques in the path connecting the two cliques in the tree. This property is called the clique intersection property, \cite{blp:94}, and trees with this property are referred to as clique trees. As a result it is possible to represent chordal graphs using clique trees. This means that in case the sparsity graph is chordal, it is possible to use algorithms for generating clique trees for chordal graphs, to extract the aforementioned tree structure in the problem. In fact this has been used for the coupled example in \eqref{eq:exampleCoupled}. Notice that the sparsity graph for this example is chordal, and the clusters marked in Figure \ref{fig:ex1b} are its cliques. Their corresponding clique tree is depicted in Figure \ref{fig:ex1c}. Also notice that in case the sparsity graph is not chordal, the same procedure can be used on its chordal embedding for extracting the tree structure. Coupled problems with a tree structure can be solved using a message-passing algorithm. Consider the following coupled convex optimization problem
\begin{align}\label{eq:coupledProblem}
\minimize_x \quad f_1(x) + f_2(x) + \dots + f_N(x),
\end{align}
where $f_i : \mathbb R^n \rightarrow \mathbb R$ for $i = 1, \dots, N$. This problem can be seen as a combination of $N$ subproblems, each of which is defined by a term in the cost function and depends only on a few elements of $x$. Let us describe the coupling structure in this problem in a similar manner as we did for the coupled SDP in \eqref{eq:coupledSDP}. That is we denote the ordered set of indices of $x$ that each subproblem $i$ depends on by $J_i$, and we denote the ordered set of indices of functions that depend on $x_j$ by $\mathcal I_j$. We can equivalently rewrite this problem as
\begin{align}\label{eq:CPS}
\minimize_{x} & \quad   \bar f_1\left(x_{_{J_1}}\right) + \dots + \bar f_N(x_{_{J_N}}),
\end{align}
where the functions $\bar f_i \ : \ \mathbb R^{|J_i|} \rightarrow \mathbb R$ are lower dimensional descriptions of $f_i$s such that $f_i(x) = \bar f_i(E_{J_i}x)$ for all $x \in \mathbb R^n$ and $i = 1, \dots, N$. Let us assume that the sparsity graph of this problem, $Q_s(V_s, \mathcal E_s)$, has an inherent tree structure with a set of cliques $\mathbf C_{Q_s} = \{ C_1, \dots, C_q \}$ and a clique tree, $T(V_t, \mathcal E_t)$. This problem can be solved distributedly using the message-passing algorithm that utilizes the clique tree as its computational graph. This means that the nodes $V_t = \{ 1, \dots, q \}$ act as computational agents that communicate or collaborate with their neighbors defined by the edge set $\mathcal E_t$. In order to describe the message-passing algorithm, we first need to assign each subproblem in \eqref{eq:CPS} to each of the agents. We can assign a subproblem or function $\bar f_i$ to an agent $j$ if $J_i \subseteq C_j$. Let us denote the set of indices of the subproblems assigned to agent $j$ by $\phi_j$. Then we can rewrite \eqref{eq:coupledProblem} as
\begin{align}\label{eq:CPSClustered}
\minimize_x \quad \sum_{i = 1}^q  F_i\left(x_{_{C_i}}\right),
\end{align}
where $F_i\left(x_{_{C_i}}\right) := \sum_{j \in \phi_i} \bar f_j\left(x_{_{J_i}}\right)$. The message-passing algorithm, much the same as dynamic programming, solves \eqref{eq:CPSClustered} by performing an upward-downward pass through the clique tree, see e.g., \cite[Sec. 4]{kho:15c}, \cite{kol:09} and references therein. Next we show how the message-passing algorithm can be used for devising distributed solvers for coupled SDPs with a tree structure.
\section{Distributed Primal-dual Interior-point Methods for Coupled SDPs}\label{sec:Dprimaldual}

Let us reconsider the coupled SDP in \eqref{eq:coupledSDP}, and assume that the sparsity graph of this problem, $Q_s(V_s, \mathcal E_s)$, has an inherent tree structure with clique set $\mathbf C_{Q_s} = \{ C_1, \dots, C_q \}$ and clique tree $T(V_t, \mathcal E_t)$. Here we propose a method that allows us to solve this problem distributedly over the clique tree. To this end, we first need to assign the constituent subproblems of~\eqref{eq:coupledSDP} to each of the agents in the tree. Firstly define $\bar C_i \subseteq \mathbb N_n$ such that $\bar C_i \times_s \bar C_i = C_i$.  Then we can assign a subproblem $i$ to agent $j$ if $J_i \subseteq \bar C_j$. As in Section \ref{sec:TreeStructureMP}, let us denote the set of indices of subproblems assigned to agent $j$ by $\phi_j$. The proposed algorithm in this section relies on primal-dual interior-point methods. As was discussed in Section \ref{sec:primaldual}, the most computationally demanding stage within the primal-dual method in Algorithm \ref{alg:PD} concerns the computation of the predictor and corrector directions. Hence the first step for devising a distributed algorithm for solving coupled SDPs is to distribute the computation of these directions, which is discussed next.

\subsection{Distributed Computation of Primal-dual Directions Using Message-passing}\label{sec:direction}

Recall that we can compute the predictor and corrector directions by solving the problem in \eqref{eq:CoupledSDPQP} for different choices of $r^{i,(k)}$. Firstly notice that the problem in \eqref{eq:CoupledSDPQP} is equivalent to the problem
\small
\begin{subequations}\label{eq:CoupledSDPQPReduced}
\begin{align}
\minimize_{\Delta x} & \quad \sum_{i=1}^N\frac{1}{2} (\Delta x^i)^T(F^{i})^{-1} U^{i} \Delta x^i + (r^{i})^T \Delta x^i\label{eq:CoupledSDPQPReduced-a}\\
\subject & \quad \mathcal Q^i \Delta x^i = r_{\textrm{primal}}^{i}, \  i = 1, \dots, N,\label{eq:CoupledSDPQPReduced-b}
\end{align}
\end{subequations}
\normalsize
with $\Delta x^i := \svec(\Delta X_{_{J_iJ_i}})$, that is achieved by eliminating the constraints in \eqref{eq:CoupledSDPQP-c}. It is then possible to compute the search directions by solving the problem in \eqref{eq:CoupledSDPQPReduced}. Particularly, by solving this problem we compute primal variables direction $\svec(\Delta X)$ and dual variables directions $\Delta \mathbf v$. Then we can construct the remaining primal and dual directions as
\begin{equation}\label{eq:remaining}
\begin{split}
\Delta \bar X^i &= \Delta X_{_{J_iJ_i}},\\
\Delta \bar v^i  &= (F^{i})^{-1}U^{i} \mathcal Q^i\Delta \bar x^i + r^{i} - (\mathcal Q^i)^T\Delta v^i.
\end{split}
\end{equation}
for $i = 1, \dots, N$. Next theorem shows that these directions in fact satisfy the system of equations in \eqref{eq:CoupledSDPAugmentedSystem}.
\begin{theorem}
The primal-dual directions computed by solving \eqref{eq:CoupledSDPQPReduced} and using \eqref{eq:remaining} satisfy the system of equations in \eqref{eq:CoupledSDPAugmentedSystem}.
\end{theorem}
\begin{IEEEproof}
Notice that any solution of \eqref{eq:CoupledSDPQPReduced} satisfies
\begin{subequations}\label{eq:thm1temp1}
\begin{align}
\bar{\mathcal E}^T \left(\mathbf F^{-1}\mathbf U \bar{\mathcal E}\Delta x - \mathcal Q^T \Delta \mathbf v \right) & = -\bar{\mathcal E}^T \mathbf r,\label{eq:thm1temp1-a}\\
\mathcal Q \bar{\mathcal E} \Delta x &= \mathbf r_{\textrm{primal}}.
\end{align}
\end{subequations}
By choosing $\Delta \bar X^i = \Delta X_{_{J_iJ_i}}$, the primal directions, $\Delta x$ and $\Delta \bar{\mathbf x}$, will satisfy the third and fourth block equations in \eqref{eq:CoupledSDPAugmentedSystem}. Furthermore, notice that by \eqref{eq:thm1temp1-a} we have that
\begin{align*}
\mathbf F^{-1}\mathbf U \Delta \bar{\mathbf x} - \mathcal Q^T \Delta \mathbf v + \mathbf r \in \mathcal N (\bar{\mathcal E}^T).
\end{align*}
So if we set
\begin{align}
\Delta \bar v^i  = (F^{i})^{-1}U^{i} \mathcal Q^i\Delta \bar x^i + r^{i} - (\mathcal Q^i)^T\Delta v^i.
\end{align}
for $i = 1, \dots, N$, not only the primal-dual iterates satisfy the first block equation in \eqref{eq:CoupledSDPAugmentedSystem}, but also we have $\bar{\mathcal E}^T \Delta \bar{\mathbf v} = 0$. This completes the proof.
\end{IEEEproof}
Consequently, we can construct the primal-dual solutions for the problem in \eqref{eq:CoupledSDPQP} by first solving the problem in \eqref{eq:CoupledSDPQPReduced} and constructing the remainder of the solution as outlined in \eqref{eq:remaining}. Notice that the coupling structure of~\eqref{eq:CoupledSDPQPReduced} is the same as that of \eqref{eq:coupledSDP}. This means that both problems have the same sparsity graph and tree representation of the coupling structure. We can equivalently rewrite \eqref{eq:CoupledSDPQPReduced} as
\begin{align}\label{eq:DirectionDecomposed}
\minimize \quad \sum_{i = 1}^N \bar f_i\left(\Delta x^i\right),
\end{align}
where $\bar f_i\left(\Delta x^i\right) := f_i\left(\Delta x^i\right) + \mathcal I_{\mathcal C_i}\left(\Delta x^i\right)$, with
\begin{align}
f_i \left(\Delta x^i\right) =  (\Delta x^i)^T (F^{i})^{-1} U^{i}\Delta x^i + (r^{i})^T \Delta x^i,
\end{align}
for $i = 1, \dots, N$, and functions $\mathcal I_{\mathcal C_i}$ for $i = 1, \dots, N$, are the indicator functions for the constraints in \eqref{eq:CoupledSDPQPReduced-b}, i.e.,
\begin{align*}
\mathcal I_{\mathcal C_i}\left(\Delta x^i\right) = \begin{cases} 0 \hspace{6mm} \mathcal Q^i\Delta x^i = r_{\textrm{primal}}^{i}  \\ \infty \hspace{4mm} \text{Otherwise} \end{cases}.
\end{align*}
This problem is in the same format as \eqref{eq:CPS}, and due to its coupling structure, can be solved distributedly using message passing, see \cite[Sec. 6.2]{kho:15c}.

So far we have described how to distribute the computation of the search directions using message passing. However, it remains to discuss how to distributedly compute the primal and dual step sizes, update the perturbation parameter and decide on terminating the algorithm. We discuss these next.

\subsection{Distributed Step-size Computation and Termination Check}\label{sec:stepsize}

The clique tree used for computing the search directions can also be used for performing the remaining computations in Algorithm~\ref{alg:PD} distributedly. Notice that the computations described in this section are different than that of presented in \cite{kho:15c}. This is because here we rely on a predictor-corrector method and we are concerned with SDPs. Let us first focus on step size computation. Similar to the message-passing algorithm, in order to compute the primal and dual step sizes we need to perform an upward-downward pass over the clique tree. We start the computation from the agents at the leaves of the tree, where every such agent first computes
\begin{subequations}
\begin{align}
\underline{\lambda}_p^i &= \minimum_{j \in \phi_i} \left( \lambda_{\minimum}\left((\bar X^{j})^{-1}\Delta \bar X^j_\textrm{pred}\right) \right),\\
\underline{\lambda}_d^i &= \minimum_{j \in \phi_i} \left( \lambda_{\minimum}\left((S^{j})^{-1}\Delta S^j_\textrm{pred}\right) \right),
\end{align}
\end{subequations}
and communicates them to its corresponding parent. Each agent $i$ that has received these quantities from their children, will then compute
\begin{subequations}\label{eq:stepsizeGen}
\begin{align}
\underline{\lambda}_p^i &= \minimum \left(\minimum_{j \in \children(i)}\left( \underline{\lambda}_p^j \right), \minimum_{j \in \phi_i}\left( \lambda_{\minimum}\left((\bar X^{j})^{-1}\Delta \bar X^j_\textrm{pred}\right)\right) \right),\\
\underline{\lambda}_d^i &= \minimum \left(\minimum_{j \in \children(i)}\left( \underline{\lambda}_d^j \right), \minimum_{j \in \phi_i}\left( \lambda_{\minimum}\left((S^{j})^{-1}\Delta S^j_\textrm{pred}\right)\right) \right),
\end{align}
\end{subequations}
and will communicate them to its parent. This procedure is then continued until we arrive at the root of the tree. At this point, the agent at the root computes the primal and dual step sizes as
\begin{align}
\alpha_p &:= \minimum \left( 1, \frac{-\tau}{\underline{\lambda}_p^r} \right), \ \alpha_d:= \minimum \left( 1, \frac{-\tau}{\underline{\lambda}_d^r} \right),
\end{align}
where $\underline{\lambda}_p^r $ and $\underline{\lambda}_d^r $ are calculated as in \eqref{eq:stepsizeGen}. These quantities are then communicated downwards through the tree until they reach the agent at the leaves. At this point, all agents will know the primal and dual step sizes. So the step sizes computation can be done by an upward-downward pass over the tree. Notice that the need for computing primal and dual step sizes also appear in Step 7 of Algorithm \ref{alg:PD}. We can use the same procedure for computing the step sizes at this step by simply replacing the predictor directions with corrector ones.

As can be seen from Algorithm \ref{alg:PD}, in order to compute the corrector directions we first need to update the parameter $\sigma$ in Step 6 of the algorithm. We can use a similar approach to perform this update distributedly over the clique tree. Let us start the computations from the leaves of the tree. Every agent $i$ at the leaves will then compute and communicate
\begin{subequations}
\begin{align}
\sigma_1^i &= \sum_{j \in \phi_i} (\bar X^{j} + \alpha_p \Delta\bar X^j_\textrm{pred})\bullet(S^{j} + \alpha_d \Delta S^j_\textrm{pred}),\\
\sigma_2^i &= \sum_{j \in \phi_i} \bar X^{j}\bullet S^{j},
\end{align}
\end{subequations}
to its corresponding parent. Then every agent $i$ that has received these quantities from its children computes and communicates
\small
\begin{subequations}\label{eq:sigmaGen}
\begin{align}
\sigma_1^i &= \sum_{j \in \children(i)} \sigma_1^j\notag \\ & \quad  \quad + \sum_{j \in \phi_i} (\bar X^{j} + \alpha_p \Delta\bar X^j_\textrm{pred})\bullet(S^{j} + \alpha_d \Delta S^j_\textrm{pred}),\\
\sigma_2^i &= \sum_{j \in \children(i)} \sigma_2^j + \sum_{j \in \phi_i} \bar X^{j}\bullet S^{j},
\end{align}
\end{subequations}
\normalsize
to its parent. This procedure is then continued until we reach the agent at the root. Then this agent also computes the quantities $\sigma_1^r$ and $\sigma_2^r$ as in \eqref{eq:sigmaGen} and calculates the update for $\sigma$ as
\begin{align}
\sigma = \left(\frac{\sigma_1^r}{\sigma_2^r}\right)^a.
\end{align}
This quantity is then communicated downwards through the tree until it reaches the leaves of the tree. Hence, at every iteration of the primal-dual method all agents will have an update of $\sigma$ after an upward-downward pass over the clique tree.

It now remains to discuss distributed computation of terms in the stopping criteria. This concerns the computation of primal and dual residuals norms together with the surrogate duality gap. These quantities can also be computed distributedly over the clique tree using an analogous approach as above. Similarly as before let us start the computations from the leaves of the tree where every such agent computes
\small
\begin{subequations}
\begin{align}
r_d^i &= \sum_{j\in \phi_i} \left\| r_{\textrm{dual}}^{j}  \right\|^2,\ r_p^i = \sum_{j\in \phi_i} \left\| r_{\textrm{primal}}^{j}  \right\|^2,\\
\mu^i &= \sum_{j \in \phi_i} \bar X^{j}\bullet S^{j},
\end{align}
\end{subequations}
\normalsize
based on the updated iterates, and communicates them to its parent. Then every agent $i$ that has received the necessary information from its children will compute
\small
\begin{subequations}
\begin{align}
r_d^i &= \sum_{j\in \children(i)} r_d^j +  \sum_{j\in \phi_i} \left\| r_{\textrm{dual}}^{j}  \right\|^2,\\
r_p^i &= \sum_{j\in \children(i)} r_p^j + \sum_{j\in \phi_i} \left\| r_{\textrm{primal}}^{j}  \right\|^2,\\
\mu^i &= \sum_{j\in \children(i)} \mu^j + \sum_{j \in \phi_i} \bar X^{j}\bullet S^{j},
\end{align}
\end{subequations}
\normalsize
based on the updated iterates, and communicates them to the respective parent. This procedure is then continued until we reach the agent at the root, which will compute the primal and dual residuals as
\small
\begin{subequations}
\begin{align}
\left\| r_{\textrm{primal}}^{(k)} \right\|^2 &= \sum_{j\in \children(r)} r_p^j + \sum_{j\in \phi_r} \left\| r_{\textrm{primal}}^{j}  \right\|^2,\\
\left\| r_{\textrm{dual}} \right\|^2 &= \sum_{j\in \children(r)} r_d^j +  \sum_{j\in \phi_r} \left\| r_{\textrm{dual}}^{j}  \right\|^2,
\end{align}
\end{subequations}
\normalsize
and the surrogate duality gap as
\small
\begin{align}
\mu &= \frac{1}{\sum_{j = 1}^N |J_j|}\left(\sum_{j\in \children(r)} \mu^j + \sum_{j \in \phi_r} \bar X^{j}\bullet S^{j}\right).
\end{align}
\normalsize
This agent will then check the stopping criteria as in Step 12 of Algorithm \ref{alg:PD}. If these criteria are satisfied, then the agent at the root will communicate the decision to terminate the algorithm downwards through the tree. Otherwise, this agent will instead communicate the surrogate duality gap. Agents will need this parameter for updating the perturbation parameter for the next iteration of the primal-dual method.

So far we have expressed how to distribute the computations in every iteration of the primal-dual method. Next we summarize the outlined distributed algorithm in this section.
\subsection{Summary of the Algorithm and Its Computational Properties}\label{sec:summary}

Let us reconsider the coupled SDP in \eqref{eq:coupledSDP}. Given such a problem and its corresponding sparsity graph, $Q_s(V_s, \mathcal E_s)$, we extract its tree structure based on clique set $C_{Q_s} = \{ \bar C_1, \dots, \bar C_q \}$. Having done so we have the computational graph for our algorithm and it is possible to assign the constituent subproblems to each of the agents using the guidelines in Section \ref{sec:TreeStructureMP} or at the beginning of Section \ref{sec:Dprimaldual}. We can now summarize our proposed distributed algorithm as below.
\begin{flushleft}
\begin{algorithmic}
\State{Given $k = 0$, $\tau \in (0,1)$, $a \in \{ 1, 2, 3 \}$, $\epsilon > 0$, $\epsilon_{\textrm{feas}} > 0$, initial iterates $X^{(0)}$, $\bar X^{i,(0)} = X^{(0)}_{_{J_iJ_i}}\succ 0$, $S^{i,{0}} \succ 0$, $v^{i,{0}}$, $\bar v^{i,{0}}$ such that $\sum_{i = 1}^N(E_{J_i} \otimes_s E_{J_i})^T \bar v^{i,(0)} = 0$, for $i = 1, \dots, N$, and $\mu = \sum_{i= 1}^N \bar X^{i,(0)}\bullet S^{i,(0)}/\left(\sum_{j = 1}^N |J_j| \right)$}
\Repeat
\For {$i = 1, \dots, q$}
\State{Agent $i$, given $\bar X^{j,(k)}$, $S^{j,(k)}$, $\bar v^{j,(k)}$, $v^{j,(k)}$, $\sigma = 0$,}
\State{and $r^{j,(k)} = r_{\textrm{dual}}^{j,(k)} - (F^{j,(k)})^{-1}r_{\textrm{cent}}^{j,(k)}$, form}
\State{$\bar f_i\left(\Delta x^i\right)$ as in \eqref{eq:DirectionDecomposed} for $j \in \phi_i$.}
\EndFor
\vspace{2mm}
\State{Perform an upward-downward pass and compute the}
\State{predictor directions using message-passing, and \eqref{eq:DeltaSi}}
\State{and \eqref{eq:remaining}.}
\vspace{2mm}
\State{Compute the primal and dual step sizes, by performing}
\State{an upward-downward pass through the tree as}
 \State{discussed in Section \ref{sec:stepsize}.}
\vspace{2mm}
\State{Update $\sigma$ by performing an upward-downward pass}
\State{through the tree as discussed in Section \ref{sec:stepsize}.}
\vspace{2mm}
\For {$i = 1, \dots, q$}
\State{Agent $i$ reforms their subproblems with $r^{j,(k)}$ as}
\State{in \eqref{eq:RCorrector} for $j \in \phi_i$.}
\EndFor
\vspace{2mm}
\State{Perform an upward-downward pass and compute the}
\State{corrector directions using message-passing, and \eqref{eq:DeltaSCorri}}
\State{and \eqref{eq:remaining}.}
\vspace{2mm}
\State{Compute the primal and dual step sizes, by performing}
\State{an upward-downward pass through the tree as}
 \State{discussed in Section \ref{sec:stepsize}.}
\vspace{2mm}
\For {$i = 1, \dots, q$}
\State{Agent $i$ updates
\small
\begin{align*}
X_{_{J_jJ_j}}^{(k+1)} &:=  X_{_{J_jJ_j}}^{(k)} + \alpha_p(\Delta X_{_{J_jJ_j}})_\textrm{corr},\\
\bar X^{j,(k+1)} &:=  \bar X^{j,(k)} + \alpha_p \Delta \bar X^j_\textrm{corr},\\
S^{j,(k+1)} &:=  S^{j,(k)} + \alpha_d\Delta S^j_\textrm{corr},\\
v^{j,(k+1)} &:=  v^{j,(k)} + \alpha_d\Delta v^j_\textrm{corr},\\
\bar v^{j,(k+1)} &:= \bar v^{j,(k)} + \alpha_d\Delta\bar v^j_\textrm{corr},
\end{align*}
\normalsize}
\State{for $j \in \phi_i$.}
\EndFor
\vspace{1mm}
\State {$k = k + 1$.}
\vspace{2mm}
\State{Evaluate $\mu$ and the termination criteria by performing}
\State{an upward-downward pass through the tree and}
\State{decide whether to terminate the algorithm.}
\vspace{2mm}
\Until{the algorithm is terminated}
\end{algorithmic}
\end{flushleft}
From the outlined algorithm, we can observe that each iteration of the primal-dual method is accomplished within six upward-downward passes through the tree. Namely, two passes for computing the predictor and corrector directions, two for computing primal and dual step sizes, one for updating $\sigma$ and one for evaluating the stopping criteria and computing the surrogate duality gap. Let the height of the tree, that is the maximum number of edges in a path from the root to a leaf, be $h$. As a result, each iteration of the primal-dual method is accomplished in $6\times 2 \times h$ steps. Furthermore, among these passes the ones required for computing the predictor and corrector directions are by far the most computationally demanding ones. This is mainly because during the upward message-passing for these passes, every agent $i$ needs to factorize a matrix, see \cite[Sec. 6.2]{kho:15c}. However, notice that at every iteration of the primal-dual method, this matrix is the same for the predictor and corrector directions computations. This means that if each agent pre-caches the factorization of this matrix during predictor directions computations, it can reuse it for corrector directions computation, see \cite[Remark~8]{kho:15c}. This significantly reduces the computational burden of the upward-downward pass for computing corrector directions. Let us assume that the primal-dual method converges within $p$ iterations. Then the major computational burden for each agent concerns the computation of $p$ factorizations of a matrix, that is commonly of comparatively small size for loosely coupled problems. This is in stark contrast to distributed algorithms that purely rely on first-order splitting methods, as at every iteration of such algorithms each agent is required to solve an SDP.
\begin{rem}\label{rem:rem2}
The algorithm presented in this section, can distributedly detect infeasibility in the sense discussed in Remark \ref{rem:rem1}, by monitoring their local primal and dual variables. In case any agent detects divergence of these variables, it can then communicate the occurrence through the tree to terminate the algorithm.
\end{rem}
Next we discuss a class of sparse SDPs, that appear in robustness analysis of large-scale interconnected uncertain systems, and we will describe how such problems can be reformulated as coupled SDPs with an inherent tree structure.

\section{Chordal Sparsity and Domain-space Decomposition}\label{sec:ChordalSDP}

In order to describe sparsity in SDPs, we first briefly discuss the use of graphs for expressing sparsity patterns of symmetric matrices.

\subsection{Sparsity and Semidefinite Matrices}

Consider a symmetric matrix $X \in \mathbb S^n$, and an undirected graph $H(V, \mathcal E)$ with $V = \{1, \dots, n\}$ and $\mathcal E = \{ (i,j) \in (V \times V)\ | \ X_{ij} \neq 0, i \neq j\}$. We refer to this graph as the sparsity pattern graph of $X$. It is also possible to use undirected graphs to describe partial symmetric matrices. A partial symmetric matrix is a symmetric matrix where only a subset of its elements are specified and the rest are free. For the symmetric matrix $X$ this structure can be expressed using $H(V, \mathcal E)$ with $V = \{ 1, \dots, n \}$ and $\mathcal E \subseteq (V\times V)$. Particularly, the edge set is such that we can express the set of indices of specified elements using $\mathbf I_s = \mathcal E \cup \{ (i,i) \ | \ i = 1, \dots, n \}$. We denote the set of partial symmetric matrices over $H(V, \mathcal E)$ by $\mathbb S^n_H$. A matrix $X \in \mathbb S_H^n$ is then said to be positive semidefinite completable if by choosing its free elements, i.e., elements with indices in $\mathbf I_f = (V \times V)\setminus \mathbf I_s$, it is possible to produce a positive semidefinite matrix. Such matrices play a central role in the upcoming discussions. Let us review a fundamental result concerning semidefinite completable matrices.
\begin{theorem}\cite[Thm. 7]{gro:84}\label{thm:NSDC}
Let $H(V, \mathcal E)$ be a chordal graph with clique set $ \mathbf C_H = \{\bar C_1, \dots, \bar C_l\}$ such that clique intersection property holds. Then $X \in \mathbb S_H^n$ is positive semidefinite completable, if and only if
\begin{align}
X_{\bar C_i\bar C_i} \succeq 0, \quad   \ i = 1, \dots, l.
\end{align}
\end{theorem}

We will next discuss how this theorem can be used for reformulating sparse SDPs.

\subsection{Domain-space Decomposition}

Consider the following inequality-form SDP
\begin{subequations}\label{eq:SparseSDP}
\begin{align}
\minimize_y & \quad c^T y\\
\subject & \quad \sum_{i = 1}^g E_{J_i}^T Q^i E_{J_i} y_i + \sum_{i = 1}^g E_{J_i}^T M^i E_{J_i} \preceq 0
\end{align}
\end{subequations}
where $y \in \mathbb R^g$, $Q^i, M^i \in \mathbb S^{|J_i|}$ and $J_i \subset \mathbb N_n$ for $i = 1, \dots, g$. Let us denote the sparsity pattern graph for the matrix $\sum_{i = 1}^gE_{J_i}^T E_{J_i}$ with $H(V, \mathcal E)$. Assume that this graph is chordal, or that we can produce a chordal embedding by adding a few edges, with clique set $ \mathbf C_H = \{\bar C_1, \dots, \bar C_l\}$. The dual problem for \eqref{eq:SparseSDP} is given as
\begin{subequations}\label{eq:SparseSDPDual2}
\begin{align}
\minimize_Z & \quad -\sum_{i = 1}^g Z_{_{J_iJ_i}} \bullet M^i\\
\subject & \quad Z_{_{J_iJ_i}} \bullet Q^i = - c_i, \quad i = 1, \dots, g,\\
& \quad Z \succeq 0.
\end{align}
\end{subequations}
We can observe that the only elements that affect the equality constraints and the cost function are the ones specified by $\mathbf I_s$. The rest are only used in the semidefinite constraint. This in turn implies that $Z \in \mathbb S^n_H$, and using Theorem \ref{thm:NSDC}, allows us to equivalently rewrite \eqref{eq:SparseSDPDual2} as
\begin{subequations}\label{eq:SparseSDPDual3}
\begin{align}
\minimize_{Z_{{\bar C_1\bar C_1}}, \dots, Z_{{\bar C_l\bar C_l}}} & \quad - \sum_{i = 1}^g Z_{_{J_iJ_i}} \bullet M^i\\
\subject & \quad Z_{_{J_iJ_i}} \bullet Q^i = - c_i, \quad i = 1, \dots, g,\label{eq:SparseSDPDual3-b}\\
& \quad Z_{_{\bar C_i \bar C_i}} \succeq 0, \quad i = 1, \dots, l.
\end{align}
\end{subequations}
This method of reformulating \eqref{eq:SparseSDPDual2} as \eqref{eq:SparseSDPDual3} is referred to as the domain-space decomposition \cite{fuk:00}, \cite{kim+koj+mev+yam10}. Notice that for every $J_i$ there exists a $\bar C_j$ such that $J_i \subseteq \bar C_j$. This is because every set $J_i$ induces a complete subgraphs on $H(V, \mathcal E)$, and hence based on the definitions of cliques, it is either a subset of a clique or a clique itself. Let us denote the set of indices of sets $J_i$ that are a subset of $\bar C_j$ by $\phi_j$. We can then group the equality constraints in \eqref{eq:SparseSDPDual3-b} and rewrite the problem in \eqref{eq:SparseSDPDual3}~as
\begin{subequations}\label{eq:SparseSDPDual4}
\begin{align}
\minimize_{Z_{{\bar C_1\bar C_1}}, \dots, Z_{{\bar C_l\bar C_l}}} & \quad -\sum_{i = 1}^g Z_{_{J_iJ_i}} \bullet M^i\\
\subject & \quad Z_{_{J_jJ_j}} \bullet Q^j = - c_j, \ j \in \phi_i, \quad i = 1, \dots, l,\\
& \quad Z_{_{\bar C_i \bar C_i}} \succeq 0, \quad i = 1, \dots, l.
\end{align}
\end{subequations}
which is in the same format as \eqref{eq:coupledSDP}. This problem comprises $l$ subproblems. Furthermore, due to its construction has a chordal sparsity graph with $q$ cliques and a clique tree that has the same structure as the clique tree for $H(V, \mathcal E)$, where instead of $\bar C_i$, the cliques are given as $\bar C_i \times_s \bar C_i$. In fact the chordality of the sparsity graph follows, since the ordering defined by the clique tree is also a perfect elimination ordering for this graph, see \cite{gol:04} for more details.
\begin{rem}\label{rem:rem3}
Notice that the discussion in this section also extends to matrices in positive semidefinite Hermitian cones, \cite{gro:84}. This means that the decomposition scheme described here, can also be used for problems with complex data matrices.
\end{rem}
Next we will discuss robustness analysis of interconnected uncertain systems and will show how the approach described here can be used for reformulating this problem as a coupled SDP.
\section{Robustness Analysis of Interconnected Uncertain Systems}\label{sec:RSA}
In this section, we discuss robustness analysis of interconnected uncertain systems using integral quadratic constraints (IQCs). We start this discussion by first reviewing the IQC analysis framework.

\subsection{Robustness Analysis using IQCs}\label{sec:IQC}
Consider the following uncertain system
%
\begin{equation}\label{eq:UncertainSystem}
p = G q, \ \ q = \Delta(p),
\end{equation}
where $G \in \mathcal{RH}_{\infty}^{m\times m}$ is the system transfer function matrix, and $\Delta: \mathbb R^m \rightarrow \mathbb R^m$ is a bounded and causal operator representing the uncertainty in the system. We can characterize the uncertainty in the system using IQCs. Particularly it is said that $\Delta$ satisfies the IQC defined by $\Pi$, i.e., $\Delta \in \IQC(\Pi)$, if
\begin{align}\label{eq:IQCT}
\int_{0}^{\infty} \begin{bmatrix} v \\ \Delta(v) \end{bmatrix}^{T} \Pi \begin{bmatrix} v \\ \Delta(v) \end{bmatrix} \, dt \geq 0, \quad \forall v \in \mathcal{L}_2^d \ ,
\end{align}
where $\Pi$ is a bounded and self-adjoint operator. This constraint can also be written in frequency domain as
\begin{align}\label{eq:IQCF}
\int_{-\infty}^{\infty} \begin{bmatrix} \widehat{v}(j\omega) \\ \widehat{\Delta(v)}(j\omega) \end{bmatrix}^{\ast} \Pi(j\omega) \begin{bmatrix} \widehat{v}(j\omega) \\ \widehat{\Delta(v)}(j\omega) \end{bmatrix} \, d\omega \geq 0,
\end{align}
where $\hat v$ and $\widehat{\Delta(v)}$ are the Fourier transforms of the signals \cite{jon:01,meg:97}. The uncertain system is then said to be robustly stable if the interconnection between $G$ and $\Delta$ remains stable for all $\Delta \in \IQC(\Pi)$. This can be established using the following theorem.
\begin{theorem}\label{thm:IQC}
The uncertain system in~\eqref{eq:UncertainSystem} is robustly stable,~if
\begin{enumerate}
\item for all $\tau \in [0,1]$ the interconnection described in \eqref{eq:UncertainSystem}, with $\tau\Delta$, is well-posed;
\item for all $\tau \in [0,1]$, $\tau \Delta \in \IQC(\Pi)$;
\item there exists $\epsilon > 0$ such that
\begin{align}\label{eq:thmIQC}
\begin{bmatrix} G(j\omega) \\ I \end{bmatrix}^{\ast} \Pi(j\omega) \begin{bmatrix} G(j\omega) \\ I \end{bmatrix} \preceq -\epsilon I,  \hspace{2mm} \forall \omega \in [0, \infty].
\end{align}
\end{enumerate}
\end{theorem}
\begin{IEEEproof}
See \cite{jon:01,meg:97}.
\end{IEEEproof}
Satisfaction of the conditions in this theorem is a sufficient condition for robustness of the uncertain system. As a result, for robustness analysis of this system it is required to find a multiplier $\Pi$ such that $\Delta \in \IQC(\Pi)$ and that it satisfies the semi-infinite LMI in \eqref{eq:thmIQC}. The condition $\Delta \in \IQC(\Pi)$ commonly imposes structural constraints on $\Pi$, and hence the analysis problem is then to find $\Pi$ with a particular structure such that it satisfies \eqref{eq:thmIQC}. It is possible to do this using either the KYP lemma, \cite{jon:01,ran:96}, or approximately using frequency-gridding, which establishes satisfaction of \eqref{eq:thmIQC} over a finite frequencies. We utilize the latter approach later as it preserves the structure in the problem. Next we describe how this framework can be used for analyzing interconnected uncertain systems.
\subsection{Robustness Analysis of Interconnected Uncertain Systems using IQCs}\label{sec:Interconnected}
An interconnected uncertain system can be viewed as a network of $N$ uncertain subsystems. We describe each of these subsystems as
\begin{equation}\label{eq:Subsystems}
\begin{split}
&p^i = G_{pq}^i q^i + G_{pw}^iw^i \\
&z^i = G_{zq}^i q^i + G_{zw}^iw^i\\
&q^i = \Delta^i(p^i),
\end{split}
\end{equation}
where $G_{pq}^i \in \mathcal{RH}_{\infty}^{d_i \times d_i}$, $G_{pw}^i \in \mathcal{RH}_{\infty}^{d_i \times m_i}$, $G_{zq}^i \in \mathcal{RH}_{\infty}^{l_i \times d_i}$, $G_{zw}^i \in \mathcal{RH}_{\infty}^{l_i \times m_i}$, and $\Delta^i:\mathbb{R}^{d_i} \to \mathbb{R}^{d_i}$. It is possible to describe the interconnection among the subsystems using a 0--1 matrix $\Gamma$ as
\begin{align}\label{eq:Interconst}
  \begin{bmatrix}
    w^1\\w^2\\ \vdots\\ w^N
  \end{bmatrix} =
\underbrace{  \begin{bmatrix}
    \Gamma_{11} & \Gamma_{12} & \cdots & \Gamma_{1N} \\
    \Gamma_{21} & \Gamma_{22} & \cdots & \Gamma_{2N} \\
    \vdots & \vdots & \ddots & \vdots \\
    \Gamma_{N1} & \Gamma_{N2} & \cdots & \Gamma_{NN}
  \end{bmatrix}}_{\Gamma}
  \begin{bmatrix}
    z^1\\z^2\\ \vdots\\ z^N
  \end{bmatrix},
\end{align}
where each $\Gamma_{ij}$ describes which components of $z^j$ is connected to which components of $w^i$. Let us define $p = (p^1, \dots, p^N)$, $q = (q^1, \dots, q^N)$, $w = (w^1, \dots, w^N)$ and $z = (z^1, \dots, z^N)$. Then we can compactly describe the entire interconnected uncertain system as
\begin{equation}\label{eq:SysInter}
\begin{split}
p& = G_{pq} q + G_{pw}w \\
z& = G_{zq} q + G_{zw}w\\
q& = \Delta(p)\\
w& = \Gamma z,
\end{split}
\end{equation}
where $G_{\star\bullet} = \diag(G_{\star\bullet}^1, \dots, G_{\star\bullet}^N)$ and $\Delta  = \diag(\Delta^1, \dots, \Delta^N)$. Let us assume that the interconnected system is nominally or internally stable, i.e., $(I-\Gamma G_{zw})^{-1} \in \mathcal{RH}_{\infty}^{\bar m \times \bar m}$ with $\bar m = \sum_{i=1}^N m_i$. It was then shown in \cite{and:13} that the system is robustly stable if there exist
\begin{align*}
\bar \Pi = \begin{bmatrix} \bar \Pi_{11} & \bar \Pi_{12} \\ \bar \Pi_{21} & \bar \Pi_{22} \end{bmatrix},
\end{align*}
with $\bar \Pi_{\star\bullet} = \diag(\Pi_{\star\bullet}^1, \dots, \Pi_{\star\bullet}^N)$ and $\Delta^i \in \IQC \left( \begin{bmatrix} \Pi^i_{11} & \Pi^i_{12} \\ \Pi^i_{21} & \Pi^i_{22}\end{bmatrix} \right)$, and a diagonal matrix $X\succ 0$ such that
\begin{multline}\label{eq:IQCInterconnected}
\begin{bmatrix} G_{pq} & G_{pw} \\ I & 0 \end{bmatrix}^{\ast}\begin{bmatrix} \bar{\Pi}_{11} & \bar{\Pi}_{12} \\ \bar{\Pi}_{21} & \bar{\Pi}_{22} \end{bmatrix}\begin{bmatrix} G_{pq} & G_{pw} \\ I & 0 \end{bmatrix} -\\
 \begin{bmatrix} -G_{zq}^{\ast}\Gamma^T\\ I
  -G_{zw}^{\ast}\Gamma^T \end{bmatrix}X\begin{bmatrix} -\Gamma
  G_{zq} & I-\Gamma G_{zw} \end{bmatrix} \preceq -\epsilon I.
\end{multline}
It is possible to rewrite this problem in the following standard form
\begin{subequations}\label{eq:SDPDual}
\begin{align}
 \find & \quad  y \label{eq:SDPDual-1} \\
 \subject & \quad  \sum_{i = 1}^{m} y_i\bar Q^i + W  \preceq 0 \label{eq:SDPDual-2}
 \end{align}
\end{subequations}
where $\bar Q^i \in \mathbb H^{\bar m + \bar d}$ for all $i = 1, \dots, m$ and $W \in \mathbb S^{\bar m + \bar d}$ with $\bar d = \sum_{i=1}^N d_i$. We can equivalently rewrite the problem in \eqref{eq:SDPDual} as below
\begin{subequations}\label{eq:SDPDual1}
\begin{align}
 \find & \quad  y \label{eq:SDPDual11} \\
 \subject & \quad  \sum_{i = 1}^{m} y_i\begin{bmatrix}\real(\bar Q^i) & -\imag(\bar Q^i) \\ \imag(\bar Q^i) & \real(\bar Q^i)\end{bmatrix} + \begin{bmatrix} W & 0 \\ 0 & W \end{bmatrix}  \preceq 0 \label{eq:SDPDual12}
 \end{align}
\end{subequations}
where all the data matrices are real, \cite{boyd:04}. In case $\Gamma$ is sparse, then this SDP is also sparse and can be written in the same format as in \eqref{eq:SparseSDP} with $c=0$. As a result we can use the approach presented in Section \ref{sec:ChordalSDP} for reformulating this problem as a coupled problem, and employ the algorithm presented in Section \ref{sec:Dprimaldual} for solving it.
\begin{rem}\label{rem:rem4}
As was discussed in Remark \ref{rem:rem3}, the decomposition can be conducted directly on \eqref{eq:IQCInterconnected} or \eqref{eq:SDPDual}. However, we here choose to reformulate the problem in \eqref{eq:SDPDual1}, with real data matrices, for ease of notation and ease of use of the algorithm described in Section \ref{sec:Dprimaldual}.
\end{rem}

Next we illustrate this approach and study the performance of the algorithm using numerical experiments.
\section{Numerical Experiments}\label{sec:numerical}

In this section we consider two examples, namely a chain of uncertain systems and an interconnected uncertain system over a so-called scale-free network. These examples are taken from \cite{and:13}.
\begin{figure}
\begin{center}
\begin{tikzpicture}[scale=0.4, inner sep=0pt]

  \tikzstyle{nodeA} = [draw, fill=none, minimum size=1.6cm]
  \tikzstyle{nodeB} = [draw, fill=none, minimum size=1.0cm]

  \draw (0,0) node [style=nodeA]
  {$G^1(s)$};
  \draw (7,0) node [style=nodeA]
  {$G^2(s)$};
  \draw (16.5,0) node [style=nodeA]
  {$G^N(s)$};

  \draw (0,4)  node [style=nodeB]
  {$\delta^1$};
  \draw (7,4)  node [style=nodeB]
  {$\delta^2$};
  \draw (16.5,4) node [style=nodeB]
  {$\delta^N$};
 
  \draw[->] (-2,1) -- (-3,1) -- node[right]{$\phantom{.}p^1$} (-3,4) -- (-1.25,4);
  \draw[->] (1.25,4) -- (3,4) -- node[left]{$q^1\phantom{.}$} (3,1) -- (2,1);
  \draw[->] (2,0) -- node[above]{\raisebox{1.5mm}{$z^1$}} (5,0);
  \draw[->] (5,-1) -- node[below]{\raisebox{-2.5mm}{$z_1^2$}} (2,-1);

  \draw[->] (5,1) -- (4,1) -- node[right]{$\phantom{.}p^2$} (4,4) -- (5.75,4);
  \draw[->] (8.25,4) -- (10,4) -- node[left]{$q^2\phantom{.}$} (10,1) -- (9,1);
  \draw[->] (9,0) -- node[above]{\raisebox{1.5mm}{$z_2^2$}} (11,0);
  \draw[->] (11,-1) --  node[below]{\raisebox{-2.5mm}{$z_1^3$}} (9,-1);

  \draw[->] (14.5,1) -- (13.5,1) -- node[right]{$\phantom{.}p^N$} (13.5,4) -- (15.25,4);
  \draw[->] (17.75,4) -- (19.5,4) -- node[left]{$q^N$} (19.5,1) -- (18.5,1);
  \draw[->] (11.75,0) -- node[above]{\raisebox{1.5mm}{$z_2^{N-1}\hspace{5mm}$}} (14.5,0);
  \draw[->] (14.5,-1) -- node[below]{\raisebox{-2.5mm}{$z^N\hspace{4mm}$}} (11.75,-1);

  \draw (11.5,-0.5) node [fill=none] {{\small  $\cdots$}};

\end{tikzpicture}
\vspace{-16pt}
\caption{A chain of $N$ uncertain subsystem.}
\label{fig:System}
\end{center}
\end{figure}
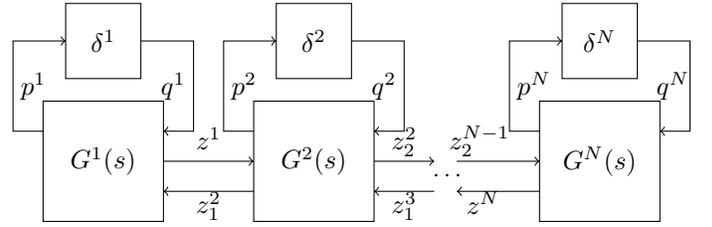
\begin{figure}[t]
\begin{center}
\includegraphics[width=9.8cm]{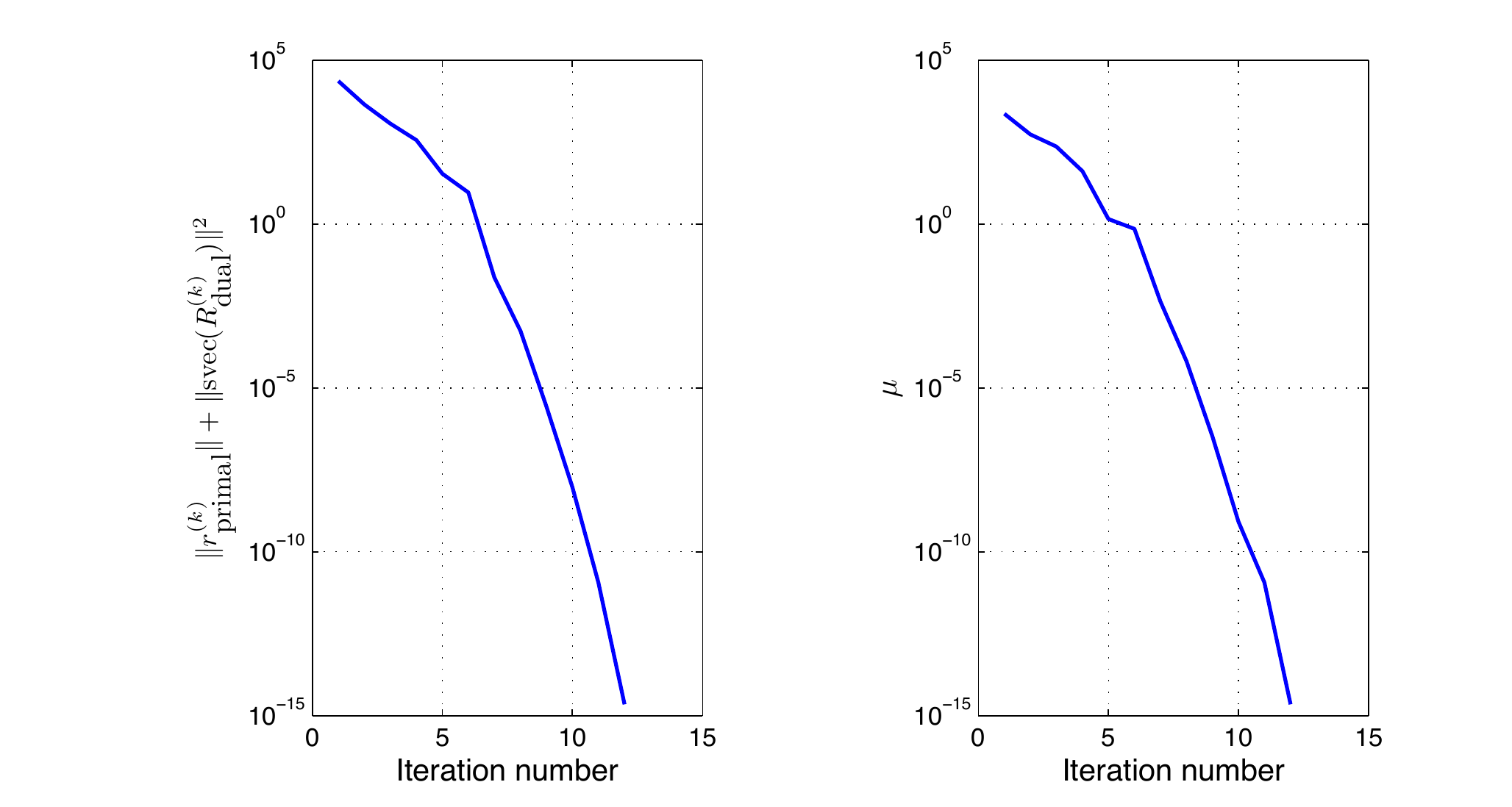}
\vspace{-16pt}
\caption{\small Convergence behavior of the algorithm for analysis of a chain of uncertain systems. The figure on the left shows the sum of primal and dual residuals and the figure on the right depicts the surrogate duality gap.\normalsize }
\label{fig:ExChain}
\end{center}
\end{figure}
Let us start with the analysis of a chain of uncertain systems, as illustrated in Figure \ref{fig:System}. As can be seen from the figure, for subsystems $1<i<N$, $z^i, w^i \in \mathbb R^2$ and for subsystems $i = 1,N$, $z^i, w^i \in \mathbb R$. The uncertainty in each subsystem $i$ is represented using $\delta^i$, which is assumed to be an unknown gain in the normalized interval $[-1, 1]$. We can hence describe the uncertainties as $\delta^i \in \IQC(\Pi^i)$ with $\Pi^i = \begin{bmatrix} r_i(j\omega) & 0 \\ 0 & -r_i(j\omega) \end{bmatrix}$, and $r_i(j\omega) \geq 0$, \cite{meg:97}. The interconnection matrix for this interconnected system is described by the nonzero blocks $\Gamma_{i,i-1} = \Gamma_{i-1,i}^T$ for $i = 2, \dots, N$, where $\Gamma_{i,i-1} = \Gamma_{i-1,i}^T = \begin{bmatrix} 0 & 1 \\ 0 & 0 \end{bmatrix}, \ i=3,\ldots,N-1$, and $ \Gamma_{21} = \Gamma_{12}^T = (1,0), \quad \Gamma_{N-1,N} =  \Gamma_{N,N-1}^T = (0,1)$. We considered the analysis problem for this system with $N = 100$ subsystems in the chain, at a single frequency $\omega = 1 \ rad/s$. We solved $10$ instances of this problem with different transfer function matrices for the subsystems. The transfer function matrices for each instance were randomly generated using the approach presented in \cite{and:13}. This guarantees that the interconnected system is robustly stable for all instances. Furthermore, for this problem the multiplier $X$ was chosen as $X = \diag(x_1, \dots, x_{2N-2})$. This resulted in a problem in the same format as in \eqref{eq:SDPDual}, with $m = 298$ and $W \in \mathbb S^{298}$.

Forming \eqref{eq:SDPDual1} for this analysis problem, resulted in an LMI with a chordal sparsity pattern, with $198$ cliques where the largest clique was of size 8. The clique tree over these cliques had a height of $99$. In order to establish chordality of the sparsity pattern graph and generate its cliques a greedy search algorithm with \emph{min degree} criterion was used, \cite{cor:01}. If we now form the problem in \eqref{eq:SparseSDPDual4}, this problem will comprise $198$ subproblems and can be solved distributedly over the clique tree. The parameters within the primal-dual method were chosen to be the same for all instances and are chosen as $a = 1$, $\tau = 0.98$, $\epsilon = \epsilon_{\textrm{feas}} = 10^{-12}$, $\bar v^{i,(0)} = 0$ and $v^{i,(0)} = 0$ for all $i = 1, \dots, N$, and $X^{(0)}$ and $S^{i,(0)}$ for $i = 1, \dots, N$ were chosen to be diagonal matrices with positive diagonal entries generated randomly with a uniform distribution in the interval $(0.1, 2)$. In the worst case the primal-dual method converged after $12$ iterations. The convergence behavior of this instance is illustrated in Figure \ref{fig:ExChain}, and as can be seen mimics that of a standard primal-dual method, i.e., convergence within 10 to 50 iterations with a quadratic convergence phase, \cite{boyd:04}. Considering the height of the tree, this algorithm then, in the worst case, converged after $6 \times 2 \times 99 \times 12 = 14256$ steps. During the run of the algorithm, each agent was required to compute a factorization $12$ times and needed to communicate with its neighbors $144$ times. The computations in the remaining steps were trivial.

\begin{figure}[t]
\begin{center}
\includegraphics[width=8.7cm,left]{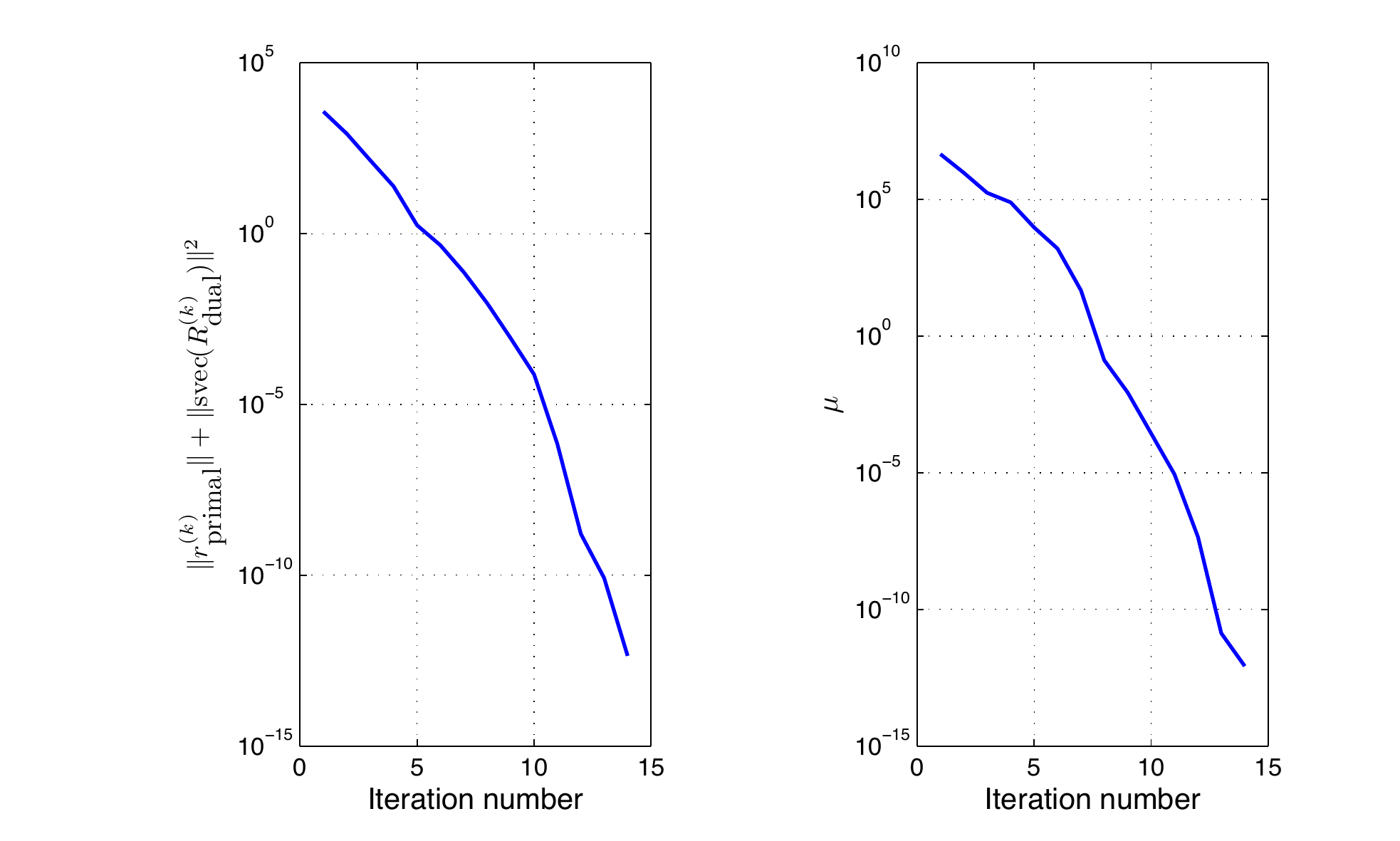}
\vspace{-18pt}   
\caption{\small Convergence behavior of the algorithm for analysis of an interconnected system over a scale-free network. The figure on the left shows the sum of primal and dual residuals and the figure on the right depicts the surrogate duality gap.\normalsize }
\label{fig:ExTree}
\end{center}
\end{figure}

We further tested the performance of the algorithm on a larger example with a more complicated interconnection description. Particularly we used the same scale-free network as in \cite[Sec. 5.2]{and:13} for describing the interconnections among the subsystems. This resulted in an extremely sparse interconnection matrix. The transfer function matrices for the subsystems were also generated using the approach presented in \cite{and:13}. Forming~\eqref{eq:SDPDual1} for this analysis problem resulted in an LMI that is sparse with $m = 1498$ and $W \in \mathbb S^{1498}$. The chordal embedding for the sparsity pattern graph of this LMI was generated by introducing $2.4\%$ fill-in, also using a greedy search algorithm, with $579$ cliques. The largest of these cliques had a size of $261$. The corresponding clique tree for this problem was of height $35$. This means that the corresponding problem in \eqref{eq:SparseSDPDual4} will comprise of $579$ subproblems and can be solved distributedly over this clique tree. We tested the performance of the proposed algorithm over 10 instances of this problem. The parameters of the primal-dual method were chosen to be the same as above. In the worst case the algorithm converged after $14$ iterations. The convergence behavior of this instance is illustrated in Figure \ref{fig:ExTree}. As a result, in the worst case, the algorithm converged after $6 \times 2 \times 35 \times 14 = 5880$ steps. During the run of the algorithm, each agent needed to compute a factorization only $14$ times and were required to communicate with its neighbors $168$ times.
\vspace{-8pt}

\section{Conclusions}\label{sec:conclusions}

In this paper we put forth a distributed algorithm for solving coupled SDPs with a tree structure. The proposed algorithm, unlike the existing ones, does not use first-order splitting methods but instead uses primal-dual interior-point methods. Particularly, this algorithm utilizes the inherent tree structure in the problem as its computational graph, and distributes the computations at each iteration of the primal-dual method among the computational agents. In order to compute the search directions at every iteration, we employ a message-passing algorithm. This enables us to compute the exact search directions in a finite number of iterations. Furthermore, we showed that this number can be computed a priori and only depends on the height of the tree. We applied the proposed algorithm for solving robustness analysis of large-scale interconnected uncertain systems, and illustrated the performance of the algorithm using numerical experiments.

As was discussed in the introduction, designing distributed algorithms are commonly conducted in two phases. Namely, a decomposition or reformulation phase and a splitting phase. In this paper, we mainly focused on the second phase of this procedure, that is design of efficient methods to distribute the computations of solving a \emph{given} coupled SDP. However, it is possible to further improve the computational and/or implementation properties of the devised algorithm, by using the available flexibilities in decomposition or reformulation phase. We will explore such possibilities as future line of research. This will mainly concern devising heuristics for clique or cluster merging to reduce the overall computational cost of the algorithm and/or to better represent the intuitive properties of the problem, such as physical structure in the problem.
\vspace{-10pt}
\ifCLASSOPTIONcaptionsoff
  \newpage
\fi



%
\bibliographystyle{plain}
\bibliography{IEEETrans}

%





\end{document}